\documentclass[11pt]{article}
\usepackage{amsmath, amsthm, amsfonts, amssymb, graphicx, fullpage, indentfirst, cite}
\theoremstyle{plain}
\newtheorem{proposition}{Proposition}
\newtheorem{lemma}{Lemma}
\newtheorem{theorem}{Theorem}
\newtheorem{corollary}{Corollary}
\theoremstyle{definition}

\newtheorem{algorithm}{Algorithm}
\theoremstyle{remark}
\newtheorem{remark}{Remark}

\newenvironment{algorithminit}[1]{\ \\{\em Initialization}: #1\begin{list}{\labelenumi}{\topsep0in\itemsep0in\parsep0in\labelwidth1in\usecounter{enumi}}}{\setcounter{enumii}{\value{enumi}}\end{list}}
\newenvironment{algorithmoper}[1]{{\em Operation}: #1\begin{list}{\labelenumi}{\topsep0in\itemsep0in\parsep0in\labelwidth1in\usecounter{enumi}\setcounter{enumi}{\value{enumii}}}}{\hfill$\blacksquare$\end{list}}

\begin{document}

\title{\Large\bf Controlled Hopwise Averaging: Bandwidth/Energy-Efficient Asynchronous Distributed Averaging for Wireless Networks\footnote{This work was supported by the National Science Foundation under grant CMMI-0900806.}}
\author{Choon Yik Tang and Jie Lu\\ School of Electrical and Computer Engineering\\ University of Oklahoma, Norman, OK 73019, USA\\ {\sf\{cytang,jie.lu-1\}@ou.edu}}
\date{\today}
\maketitle

\begin{abstract}
This paper addresses the problem of averaging numbers across a wireless network from an important, but largely neglected, viewpoint: {\em bandwidth/energy efficiency}. We show that existing distributed averaging schemes have several drawbacks and are inefficient, producing networked dynamical systems that evolve with wasteful communications. Motivated by this, we develop {\em Controlled Hopwise Averaging} (CHA), a distributed asynchronous algorithm that attempts to ``make the most'' out of each iteration by fully exploiting the broadcast nature of wireless medium and enabling control of when to initiate an iteration. We show that CHA admits a common quadratic Lyapunov function for analysis, derive bounds on its exponential convergence rate, and show that they outperform the convergence rate of Pairwise Averaging for some common graphs. We also introduce a new way to apply Lyapunov stability theory, using the Lyapunov function to perform greedy, decentralized, feedback iteration control. Finally, through extensive simulation on random geometric graphs, we show that CHA is substantially more efficient than several existing schemes, requiring far fewer transmissions to complete an averaging task.
\end{abstract}

\section{Introduction}\label{sec:intro}

Averaging numbers across a network is a need that arises in many applications of mobile ad hoc networks and wireless sensor networks. In order to collaboratively accomplish a task, nodes often have to compute the network-wide average of their individual observations. For examples, by averaging their individual throughputs, an ad hoc network of computers can assess how well the network, as a whole, is performing, and by averaging their humidity measurements, a wireless network of sensing agents can cooperatively detect the occurrence of local, deviation-from-average anomalies. Therefore, methods that enable such computation are of notable interest. Moreover, for performance reasons, it is desirable that the methods developed be robust, scalable, and efficient.

In principle, computation of network-wide averages may be accomplished via {\em flooding}, whereby every node floods the network with its observation, as well as {\em centralized computation}, whereby a central node uses an overlay tree to collect all the node observations, calculate their average, and send it back to every node. These two methods, unfortunately, have serious limitations: flooding is extremely bandwidth and energy inefficient because it propagates redundant information across the network, ignoring the fact that the ultimate goal is to simply determine the average. Centralized computation, on the other hand, is vulnerable to node mobility, node membership changes, and single-point failures, making it necessary to frequently maintain the overlay tree and occasionally start over with a new central node, both of which are rather costly to implement.

The limitations of flooding and centralized computation have motivated the search for {\em distributed averaging} algorithms that require neither flooding of node observations, nor construction of overlay trees and routing tables, to execute. To date, numerous such algorithms have been developed in continuous-time \cite{Olfati-Saber04, Cortes06, Tahbaz-Salehi07b} as well as in discrete-time for both synchronous \cite{Kempe03, Olfati-Saber04, XiaoL04, Scherber05, Kingston06, Olshevsky06, Tahbaz-Salehi07b, XiaoL07, Fagnani08, ZhuM08} and asynchronous \cite{Tsitsiklis84, Jelasity04, Montresor04, Boyd06, CaoM06, ChenJY06, Moallemi06, Mehyar07, Fagnani08} models. The closely related topic of {\em distributed consensus}, where nodes seek to achieve an arbitrary network-wide consensus on their individual opinions, has also been extensively studied; see \cite{Bertsekas89, Lynch96} for early treatments, \cite{Jadbabaie03, Olfati-Saber04, Hatano05, Moreau05, RenW05, FangL06, Sundaram07, Fagnani08, Olshevsky08, Tahbaz-Salehi10} for more recent work, and \cite{Olfati-Saber07} for a survey.

Although the current literature offers a rich collection of distributed averaging schemes along with in-depth analysis of their behaviors, their efficacy from a {\em bandwidth/energy efficiency} standpoint has not been examined. This paper is devoted to studying the distributed averaging problem from this standpoint. Its contributions are as follows: we first show that the existing schemes---regardless of whether they are developed in continuous- or discrete-time, for synchronous or asynchronous models---have a few deficiencies and are inefficient, producing networked dynamical systems that evolve with wasteful communications. To address these issues, we develop {\em Random Hopwise Averaging} (RHA), an asynchronous distributed averaging algorithm with several positive features, including a novel one among the asynchronous schemes: an ability to fully exploit the broadcast nature of wireless medium, so that no overheard information is ever wastefully discarded. We show that RHA admits a common quadratic Lyapunov function, is almost surely asymptotically convergent, and eliminates all but one of the deficiencies facing the existing schemes.

To tackle the remaining deficiency, on lack of control, we introduce the concept of {\em feedback iteration control}, whereby individual nodes use feedback to control when to initiate an iteration. Although simple and intuitive, this concept, somewhat surprisingly, has not been explored in the literature on distributed averaging \cite{Tsitsiklis84, Kempe03, Jelasity04, Montresor04, Olfati-Saber04, XiaoL04, Scherber05, Boyd06, CaoM06, ChenJY06, Cortes06, Kingston06, Moallemi06, Olshevsky06, Mehyar07, Tahbaz-Salehi07b, XiaoL07, Fagnani08, ZhuM08} and distributed consensus \cite{Bertsekas89, Lynch96, Jadbabaie03, Olfati-Saber04, Hatano05, Moreau05, RenW05, FangL06, Olfati-Saber07, Sundaram07, Fagnani08, Olshevsky08, Tahbaz-Salehi10}. We show that RHA, along with the common quadratic Lyapunov function, exhibits features that enable a greedy, decentralized approach to feedback iteration control, which leads to bandwidth/energy-efficient iterations at zero feedback cost. Based on this approach, we present two modified versions of RHA: an ideal version referred to as {\em Ideal Controlled Hopwise Averaging} (ICHA), and a practical one referred to simply as {\em Controlled Hopwise Averaging} (CHA). We show that ICHA yields a networked dynamical system with state-dependent switching, derive deterministic bounds on its exponential convergence rate for general and specific graphs, and show that the bounds are better than the stochastic convergence rate of Pairwise Averaging \cite{Tsitsiklis84, Fagnani08} for path, cycle, and complete graphs. We also show that CHA is able to closely mimic the behavior of ICHA, achieving the same bounds on its convergence rate. Finally, via extensive simulation on random geometric graphs, we demonstrate that CHA is substantially more bandwidth/energy efficient than Pairwise Averaging \cite{Tsitsiklis84}, Consensus Propagation \cite{Moallemi06}, Algorithm A2 of \cite{Mehyar07}, and Distributed Random Grouping \cite{ChenJY06}, requiring far fewer transmissions to complete an averaging task. In particular, CHA is twice more efficient than the most efficient existing scheme when the network is sparsely connected.

The outline of this paper is as follows: Section~\ref{sec:probform} formulates the distributed averaging problem. Section~\ref{sec:defiexissche} describes the deficiencies of the existing schemes. Sections~\ref{sec:RHA} and~\ref{sec:CHA} develop RHA and CHA and characterize their convergence properties. In Section~\ref{sec:perfcomp}, their comparison with several existing schemes is carried out. Finally, Section~\ref{sec:concl} concludes the paper. The proofs of the main results are included in Appendix~\ref{sec:app}.

\section{Problem Formulation}\label{sec:probform}

Consider a multi-hop wireless network consisting of $N\ge2$ nodes, connected by $L$ bidirectional links in a fixed topology. The network is modeled as a connected, undirected graph $\mathcal{G}=(\mathcal{V},\mathcal{E})$, where $\mathcal{V}=\{1,2,\ldots,N\}$ represents the set of $N$ nodes (vertices) and $\mathcal{E}\subset\{\{i,j\}:i,j\in\mathcal{V},i\neq j\}$ represents the set of $L$ links (edges). Any two nodes $i,j\in\mathcal{V}$ are one-hop neighbors and can communicate if and only if $\{i,j\}\in\mathcal{E}$. The set of one-hop neighbors of each node $i\in\mathcal{V}$ is denoted as $\mathcal{N}_i=\{j\in\mathcal{V}:\{i,j\}\in\mathcal{E}\}$, and the communications are assumed to be delay- and error-free, with no quantization. Each node $i\in\mathcal{V}$ observes a scalar $y_i\in\mathbb{R}$, and all the $N$ nodes wish to determine the network-wide average $x^*\in\mathbb{R}$ of their individual observations, given by
\begin{align}
x^*=\frac{1}{N}\sum_{i\in\mathcal{V}}y_i.\label{eq:x*=1Nsumy}
\end{align}

Given the above model, the problem addressed in this paper is how to construct a distributed averaging algorithm---continuous- or discrete-time, synchronous or otherwise---with which each node $i\in\mathcal{V}$ repeatedly communicates with its one-hop neighbors, iteratively updates its estimate $\hat{x}_i\in\mathbb{R}$ of the unknown average $x^*$ in \eqref{eq:x*=1Nsumy}, and asymptotically drives $\hat{x}_i$ to $x^*$---all while consuming bandwidth and energy efficiently.

The bandwidth/energy efficiency of an algorithm is measured by {\em the number of real-number transmissions it needs to drive all the $\hat{x}_i$'s to a sufficiently small neighborhood of $x^*$}, essentially completing the averaging task. This quantity is a natural measure of efficiency because the smaller it is, the lesser bandwidth is occupied, the lesser energy is expended for communications, and the faster an averaging task may be completed. These, in turn, imply more bandwidth and time for other tasks, smaller probability of collision, longer lifetime for battery-powered nodes, and possible earlier return to sleep mode, all of which are desirable. The quantity also allows algorithms with different numbers of real-number transmissions per iteration to be fairly compared. Although, in networking, every message inevitably contains overhead (e.g., transmitter/receiver IDs and message type), we exclude such overhead when measuring efficiency since it is not inherent to an algorithm, may be reduced by piggybacking messages, and becomes negligible when averaging long vectors.

\section{Deficiencies of Existing Schemes}\label{sec:defiexissche}

As was pointed out in Section~\ref{sec:intro}, the current literature offers a variety of distributed averaging schemes for solving the problem formulated in Section~\ref{sec:probform}. Unfortunately, as is explained below, they suffer from a number of deficiencies, especially a lack of bandwidth/energy efficiency, by producing networked dynamical systems that evolve with wasteful real-number transmissions.

The {\em continuous-time} algorithms in \cite{Olfati-Saber04, Cortes06, Tahbaz-Salehi07b} have the following deficiency:
\begin{enumerate}
\renewcommand{\theenumi}{D\arabic{enumi}}\itemsep-\parsep
\item {\em Costly discretization}: As immensely inefficient as flooding is, the continuous-time algorithms in \cite{Olfati-Saber04, Cortes06, Tahbaz-Salehi07b} may be more so: flooding only requires $N^2$ real-number transmissions for all the $N$ nodes to exactly determine the average $x^*$ (since it takes $N$ real-number transmissions for each node $i\in\mathcal{V}$ to flood the network with its $y_i$), whereas these algorithms may need far more than that to essentially complete an averaging task. For instance, the algorithm in \cite{Olfati-Saber04} updates the estimates $\hat{x}_i$'s of $x^*$ according to the differential equation
\begin{align}
\frac{d\hat{x}_i(t)}{dt}=\sum_{j\in\mathcal{N}_i}(\hat{x}_j(t)-\hat{x}_i(t)),\quad\forall i\in\mathcal{V}.\label{eq:dxhdt=sumxhxh}
\end{align}
To realize \eqref{eq:dxhdt=sumxhxh}, each node $i\in\mathcal{V}$ has to continuously monitor the $\hat{x}_j(t)$ of every one-hop neighbor $j\in\mathcal{N}_i$. If this can be done without wireless communications (e.g., by direct sensing), then the bandwidth/energy efficiency issue is moot. If wireless communications must be employed, then \eqref{eq:dxhdt=sumxhxh} has to be discretized, either exactly via a zero-order hold, i.e.,
\begin{align}
\hat{x}_i((k+1)T)=\sum_{j\in\mathcal{V}}h_{ij}\hat{x}_j(kT),\quad\forall i\in\mathcal{V},\label{eq:xh=sumhxh}
\end{align}
or approximately via numerical techniques such as the Euler forward difference method, i.e.,
\begin{align}
\frac{\hat{x}_i((k+1)T)-\hat{x}_i(kT)}{T}=\sum_{j\in\mathcal{N}_i}(\hat{x}_j(kT)-\hat{x}_i(kT)),\quad\forall i\in\mathcal{V},\label{eq:xhxhT=sumxhxh}
\end{align}
where each $h_{ij}\in\mathbb{R}$ is the $ij$-entry of $e^{-\mathbf{L}T}$, $\mathbf{L}\in\mathbb{R}^{N\times N}$ is the Laplacian matrix of the graph $\mathcal{G}$ that governs the dynamics \eqref{eq:dxhdt=sumxhxh}, and $T>0$ is the sampling period. Regardless of \eqref{eq:xh=sumhxh} or \eqref{eq:xhxhT=sumxhxh}, they may be far more costly to realize than flooding: with \eqref{eq:xh=sumhxh}, $N^2$ real-number transmissions are already needed per iteration (since, in general, $h_{ij}\neq0$ $\forall i,j\in\mathcal{V}$, so that each node $i\in\mathcal{V}$ has to flood the network with its $\hat{x}_i(kT)$, for every $k$). In contrast, with \eqref{eq:xhxhT=sumxhxh}, only $N$ real-number transmissions are needed per iteration (since each node $i\in\mathcal{V}$ only has to wirelessly transmit its $\hat{x}_i(kT)$ once, to every one-hop neighbor $j\in\mathcal{N}_i$, for every $k$). However, the number of iterations, needed for all the $\hat{x}_i(kT)$'s to converge to an acceptable neighborhood of $x^*$, may be very large, since $T$ must be sufficiently small for \eqref{eq:xhxhT=sumxhxh} to be stable. If this number exceeds $N$---which is possible and likely so with a conservatively small $T$---then \eqref{eq:xhxhT=sumxhxh} would be worse than flooding (flooding is, of course, more storage and bookkeeping intensive).\label{enu:costdisc}
\end{enumerate}

The {\em discrete-time synchronous} algorithms in \cite{Kempe03, Olfati-Saber04, XiaoL04, Scherber05, Kingston06, Olshevsky06, Tahbaz-Salehi07b, XiaoL07, Fagnani08, ZhuM08} have the following deficiencies:
\begin{enumerate}
\renewcommand{\theenumi}{D\arabic{enumi}}\itemsep-\parsep\setcounter{enumi}{1}
\item {\em Clock synchronization}: The discrete-time synchronous algorithms in \cite{Kempe03, Olfati-Saber04, XiaoL04, Scherber05, Kingston06, Olshevsky06, Tahbaz-Salehi07b, XiaoL07, Fagnani08, ZhuM08} require all the $N$ nodes to always have the same clock to operate. Although techniques for reducing clock synchronization errors are available, it is still desirable that this requirement can be removed.\label{enu:clocsync}
\item {\em Forced transmissions}: The algorithms in \cite{Olfati-Saber04, XiaoL04, Scherber05, Kingston06, Olshevsky06, Tahbaz-Salehi07b, XiaoL07, Fagnani08} update the estimates $\hat{x}_i$'s of $x^*$ according to the difference equation
\begin{align}
\hat{x}_i(k+1)=w_{ii}(k)\hat{x}_i(k)+\sum_{j\in\mathcal{N}_i}w_{ij}(k)\hat{x}_j(k),\quad\forall i\in\mathcal{V},\label{eq:xh=wxhsumwxh}
\end{align}
where each $w_{ij}(k)\in\mathbb{R}$ is a weighting factor that is typically constant. The $w_{ij}(k)$'s may be specified in several ways, including choosing them to maximize the convergence rate \cite{XiaoL04} or minimize the mean-square deviation \cite{XiaoL07}. However, no matter how the $w_{ij}(k)$'s are chosen, these algorithms are bandwidth/energy inefficient because the underlying update rule \eqref{eq:xh=wxhsumwxh} simply forces every node $i\in\mathcal{V}$ at each iteration $k$ to transmit its $\hat{x}_i(k)$ to its one-hop neighbors, irrespective of whether such transmissions are worthy. It is possible, for example, that the $\hat{x}_i(k)$'s of a cluster of nearby nodes are almost equal, so that their $\hat{x}_i(k+1)$'s, being convex combinations of their $\hat{x}_i(k)$'s, are also almost equal, causing their transmissions to be unworthy. The fact that $N$ real-number transmissions are needed per iteration also implies that \eqref{eq:xh=wxhsumwxh} must drive all the $\hat{x}_i(k)$'s to an acceptable neighborhood of $x^*$ within at most $N$ iterations, in order to just outperform flooding.\label{enu:forctran}
\item {\em Computing intermediate quantities}: The scheme in \cite{Olshevsky06} uses two parallel runs of a consensus algorithm to obtain two consensus values and defines each $\hat{x}_i(k)$ as the ratio of these two values. While possible, this scheme is likely inefficient because it attempts to compute two {\em intermediate} quantities, as opposed to computing $x^*$ directly.\label{enu:compintequan}
\end{enumerate}

The {\em discrete-time asynchronous} algorithms in \cite{Tsitsiklis84, Jelasity04, Montresor04, Boyd06, CaoM06, ChenJY06, Moallemi06, Mehyar07, Fagnani08} have the following deficiencies:
\begin{enumerate}
\renewcommand{\theenumi}{D\arabic{enumi}}\itemsep-\parsep\setcounter{enumi}{4}
\item {\em Wasted receptions}: Each iteration of Pairwise Averaging \cite{Tsitsiklis84}, Anti-Entropy Aggregation \cite{Jelasity04, Montresor04}, Randomized Gossip Algorithm \cite{Boyd06}, and Accelerated Gossip Algorithm \cite{CaoM06} involves a pair of nodes transmitting to each other their state variables. Due to the broadcast nature of wireless medium, their transmissions are overheard by unintended nearby nodes, who would immediately discard this ``free'' information, instead of using it to possibly speed up convergence, enhancing bandwidth/energy efficiency. Hence, these algorithms result in wasted receptions. The same can be said about Consensus Propagation \cite{Moallemi06} and Algorithm A2 of \cite{Mehyar07}, although they do not assume pairwise exchanges. It can also be said about Distributed Random Grouping \cite{ChenJY06}, which only slightly exploits such broadcast nature: the leader of a group does, but the members, who contribute the majority of the transmissions, do not.\label{enu:wastrece}
\item {\em Overlapping iterations}: Pairwise Averaging \cite{Tsitsiklis84}, Anti-Entropy Aggregation \cite{Jelasity04, Montresor04}, Randomized Gossip Algorithm \cite{Boyd06}, Accelerated Gossip Algorithm \cite{CaoM06}, and Distributed Random Grouping \cite{ChenJY06} require sequential transmissions from multiple nodes to execute an iteration. This suggests that before an iteration completes, the nodes involved may be asked to participate in other iterations initiated by those unaware of the ongoing iteration. Thus, these algorithms are prone to overlapping iterations and, therefore, to deadlock situations \cite{Mehyar07}. It is noted that this practical issue is naturally avoided by Consensus Propagation \cite{Moallemi06} and explicitly handled by Algorithms A1 and A2 of \cite{Mehyar07}.\label{enu:overiter}
\item {\em Uncontrolled iterations}: The discrete-time asynchronous algorithms in \cite{Tsitsiklis84, Jelasity04, Montresor04, Boyd06, CaoM06, ChenJY06, Moallemi06, Mehyar07} do not let individual nodes use information available to them during runtime (e.g., history of the state variables they locally maintain) to control when to initiate an iteration and who to include in the iteration. Indeed, Pairwise Averaging \cite{Tsitsiklis84}, Anti-Entropy Aggregation \cite{Jelasity04, Montresor04}, Accelerated Gossip Algorithm \cite{CaoM06}, Consensus Propagation \cite{Moallemi06}, and Algorithm A2 of \cite{Mehyar07} focus mostly on how nodes would update their state variables during an iteration, saying little about how they could use such information to control the iterations. Randomized Gossip Algorithm \cite{Boyd06} and Distributed Random Grouping \cite{ChenJY06}, on the other hand, let nodes randomly initiate an iteration according to some probabilities. Although these probabilities may be optimized \cite{Boyd06, ChenJY06}, the optimization is carried out {\em a priori}, dependent only on the graph $\mathcal{G}$ and independent of the nodes' state variables during runtime. Consequently, wasteful iterations may occur, despite the optimality. For instance, suppose Randomized Gossip Algorithm \cite{Boyd06} is utilized, and a pair of adjacent nodes $i,j\in\mathcal{V}$ have just finished gossiping with each other, so that $\hat{x}_i$ and $\hat{x}_j$ are equal. Since the optimal probabilities are generally nonzero, nodes $i$ and $j$ may gossip with each other again before any of them gossips with someone else, causing $\hat{x}_i$ and $\hat{x}_j$ to remain unchanged, wasting that particular gossip. Similarly, suppose Distributed Random Grouping \cite{ChenJY06} is employed, and a node $i\in\mathcal{V}$ has just finished leading an iteration, so that $\hat{x}_i$ and $\hat{x}_j$ $\forall j\in\mathcal{N}_i$ are equal. Due again to nonzero probabilities, node $i$ may lead another iteration before any of its one- or two-hop neighbors leads an iteration, causing $\hat{x}_i$ and $\hat{x}_j$ $\forall j\in\mathcal{N}_i$ to stay the same, wasting that particular iteration. These examples suggest that not letting nodes control the iterations is detrimental to bandwidth/energy efficiency and, conceivably, letting them do so may cut down on wasteful iterations, improving efficiency.\label{enu:uncoiter}
\item {\em Steady-state errors}: Consensus Propagation \cite{Moallemi06} ensures that all the $\hat{x}_i$'s asymptotically converge to the same steady-state value. However, this value is, in general, not equal to $x^*$ (see Figure~\ref{fig:xhat_nrt} of Section~\ref{sec:perfcomp} for an illustration). Although the error can be made arbitrarily small, it comes at the expense of increasingly slow convergence \cite{Moallemi06}, which is undesirable.\label{enu:steastaterro}
\item {\em Lack of convergence guarantees}: Accelerated Gossip Algorithm \cite{CaoM06}, developed based on the power method in numerical analysis, is shown by simulation to have the potential of speeding up the convergence of Randomized Gossip Algorithm \cite{Boyd06} by a factor of $10$. Furthermore, whenever all the $\hat{x}_i$'s converge, they must converge to $x^*$. However, it was not established in \cite{CaoM06} that they would always converge.\label{enu:lackconvguar}
\end{enumerate}

\section{Random Hopwise Averaging}\label{sec:RHA}

Deficiencies~\ref{enu:costdisc}--\ref{enu:lackconvguar} facing the existing distributed averaging schemes raise a question: {\em is it possible to develop an algorithm, which does not at all suffer from these deficiencies?} In this section, we construct an algorithm that simultaneously eliminates all but issue~\ref{enu:uncoiter} with uncontrolled iterations. In the next section, we will modify the algorithm to address this issue.

To circumvent the costly discretization issue~\ref{enu:costdisc} facing the existing continuous-time algorithms and the clock synchronization and forced transmissions issues~\ref{enu:clocsync} and~\ref{enu:forctran} facing the existing discrete-time synchronous algorithms, the algorithm we construct must be {\em asynchronous}, regardless of whether the nodes have access to the same global clock. To avoid issue~\ref{enu:overiter} with overlapping iterations, each iteration of this algorithm must involve only a {\em single} node sending a {\em single} message to its one-hop neighbors, without needing them to reply. To tackle issue~\ref{enu:wastrece} with wasted receptions, all the neighbors, upon hearing the same message, have to ``meaningfully'' incorporate it into updating their state variables, rather than simply discarding it. To overcome issues~\ref{enu:steastaterro} and~\ref{enu:lackconvguar} with steady-state errors and convergence guarantees, the algorithm must be asymptotically convergent to the correct average. Finally, to eliminate~\ref{enu:compintequan}, it has to avoid computing intermediate quantities.

To develop an algorithm having the aforementioned properties, consider a networked dynamical system, defined on the graph $\mathcal{G}=(\mathcal{V},\mathcal{E})$ as follows: associated with each link $\{i,j\}\in\mathcal{E}$ are a parameter $c_{\{i,j\}}>0$ and a state variable $x_{\{i,j\}}\in\mathbb{R}$ of the system. In addition, associated with each node $i\in\mathcal{V}$ is an output variable $\hat{x}_i\in\mathbb{R}$, which represents its estimate of the unknown average $x^*$ in \eqref{eq:x*=1Nsumy}. Since the graph $\mathcal{G}$ has $L$ links and $N$ nodes, the system has $L$ parameters $c_{\{i,j\}}$'s, $L$ state variables $x_{\{i,j\}}$'s, and $N$ output variables $\hat{x}_i$'s. To describe the system dynamics, let $x_{\{i,j\}}(0)$ and $\hat{x}_i(0)$ represent the initial values of $x_{\{i,j\}}$ and $\hat{x}_i$, and $x_{\{i,j\}}(k)$ and $\hat{x}_i(k)$ their values upon completing each iteration $k\in\mathbb{P}$, where $\mathbb{P}$ denotes the set of positive integers. With these notations, the state and output equations governing the system dynamics may be stated as
\begin{align}
x_{\{i,j\}}(k)&=\begin{cases}\dfrac{\sum_{\ell\in\mathcal{N}_{u(k)}}c_{\{u(k),\ell\}}x_{\{u(k),\ell\}}(k-1)}{\sum_{\ell\in\mathcal{N}_{u(k)}}c_{\{u(k),\ell\}}}, & \text{if $u(k)\in\{i,j\}$},\\ x_{\{i,j\}}(k-1), & \text{otherwise},\end{cases}\quad\forall k\in\mathbb{P},\;\forall\{i,j\}\in\mathcal{E},\label{eq:x=sumcxsumcuijx}\displaybreak[0]\\
\hat{x}_i(k)&=\frac{\sum_{j\in\mathcal{N}_i}c_{\{i,j\}}x_{\{i,j\}}(k)}{\sum_{j\in\mathcal{N}_i}c_{\{i,j\}}},\quad\forall k\in\mathbb{N},\;\forall i\in\mathcal{V},\label{eq:xh=sumcxsumc}
\end{align}
where $u(k)\in\mathcal{V}$ is a variable to be interpreted shortly and $\mathbb{N}$ denotes the set of nonnegative integers. Equation \eqref{eq:xh=sumcxsumc} says that the output variable associated with each node is a convex combination of the state variables associated with links incident to the node. Equation \eqref{eq:x=sumcxsumcuijx} says that at each iteration $k\in\mathbb{P}$, the state variables associated with links incident to node $u(k)$ are set equal to the same convex combination of their previous values. Equation \eqref{eq:x=sumcxsumcuijx} also implies that the system is a linear switched system, since \eqref{eq:x=sumcxsumcuijx} may be written as
\begin{align}
\mathbf{x}(k)=\mathbf{A}_{u(k)}\mathbf{x}(k-1),\quad\forall k\in\mathbb{P},\label{eq:x=Aux}
\end{align}
where $\mathbf{x}(k)\in\mathbb{R}^L$ is the state vector obtained by stacking the $L$ $x_{\{i,j\}}(k)$'s, $\mathbf{A}_{u(k)}\in\mathbb{R}^{L\times L}$ is a time-varying matrix taking one of $N$ possible values $\mathbf{A}_1,\mathbf{A}_2,\ldots,\mathbf{A}_N$ depending on $u(k)$, and each $\mathbf{A}_i\in\mathbb{R}^{L\times L}$ is a row stochastic matrix whose entries depend on the $c_{\{i,j\}}$'s. Hence, the sequence $(u(k))_{k=1}^\infty$ fully dictates how the asynchronous iteration \eqref{eq:x=sumcxsumcuijx} takes place, or equivalently, how the system \eqref{eq:x=Aux} switches. Throughout this section, we assume that $(u(k))_{k=1}^\infty$ is an independent and identically distributed random sequence with a uniform distribution, i.e.,
\begin{align}
\operatorname{P}\{u(k)=i\}=\frac{1}{N},\quad\forall k\in\mathbb{P},\;\forall i\in\mathcal{V}.\label{eq:Pu=i=1N}
\end{align}

\begin{remark}\label{rem:alterandequi}
Clearly, alternatives to letting $(u(k))_{k=1}^\infty$ be random and equiprobable are possible, and perhaps beneficial. We will explore such alternatives in Section~\ref{sec:CHA}, when we discuss control.\hfill$\blacksquare$
\end{remark}

For the system \eqref{eq:x=sumcxsumcuijx}, \eqref{eq:xh=sumcxsumc}, \eqref{eq:Pu=i=1N} to solve the distributed averaging problem, the $\hat{x}_i(k)$'s must asymptotically approach $x^*$ of \eqref{eq:x*=1Nsumy}, i.e.,
\begin{align}
\lim_{k\rightarrow\infty}\hat{x}_i(k)=x^*,\quad\forall i\in\mathcal{V}.\label{eq:limxh=x*}
\end{align}
Due to \eqref{eq:xh=sumcxsumc}, condition \eqref{eq:limxh=x*} is met if the $x_{\{i,j\}}(k)$'s satisfy
\begin{align}
\lim_{k\rightarrow\infty}x_{\{i,j\}}(k)=x^*,\quad\forall\{i,j\}\in\mathcal{E}.\label{eq:limx=x*}
\end{align}
To ensure \eqref{eq:limx=x*}, the parameters $c_{\{i,j\}}$'s and initial states $x_{\{i,j\}}(0)$'s must satisfy a condition. To derive the condition, observe from \eqref{eq:x=sumcxsumcuijx} that no matter what $u(k)$ is, the expression $\sum_{\{i,j\}\in\mathcal{E}}c_{\{i,j\}}x_{\{i,j\}}(k)$ is conserved after every iteration $k\in\mathbb{P}$, i.e.,
\begin{align}
\sum_{\{i,j\}\in\mathcal{E}}c_{\{i,j\}}x_{\{i,j\}}(k)=\sum_{\{i,j\}\in\mathcal{E}}c_{\{i,j\}}x_{\{i,j\}}(k-1),\quad\forall k\in\mathbb{P}.\label{eq:sumcx=sumcx}
\end{align}
Therefore, as it follows from \eqref{eq:sumcx=sumcx} and \eqref{eq:x*=1Nsumy}, \eqref{eq:limx=x*} holds only if the $c_{\{i,j\}}$'s and $x_{\{i,j\}}(0)$'s satisfy
\begin{align}
\frac{\sum_{\{i,j\}\in\mathcal{E}}c_{\{i,j\}}x_{\{i,j\}}(0)}{\sum_{\{i,j\}\in\mathcal{E}}c_{\{i,j\}}}=\frac{\sum_{i\in\mathcal{V}}y_i}{N}.\label{eq:sumcxsumc=sumyN}
\end{align}
To achieve \eqref{eq:sumcxsumc=sumyN}, notice that the expressions $\sum_{\{i,j\}\in\mathcal{E}}c_{\{i,j\}}$ and $\sum_{\{i,j\}\in\mathcal{E}}c_{\{i,j\}}x_{\{i,j\}}(0)$ each has $L$ terms, of which $|\mathcal{N}_i|$ terms are associated with links incident to node $i$, for every $i\in\mathcal{V}$, where $|\cdot|$ denotes the cardinality of a set. Hence, by letting each node $i\in\mathcal{V}$ evenly distribute the number $1$ to the $|\mathcal{N}_i|$ terms in $\sum_{\{i,j\}\in\mathcal{E}}c_{\{i,j\}}$, i.e.,
\begin{align}
c_{\{i,j\}}=\frac{1}{|\mathcal{N}_i|}+\frac{1}{|\mathcal{N}_j|},\quad\forall\{i,j\}\in\mathcal{E},\label{eq:c=1|N|1|N|}
\end{align}
we get $\sum_{\{i,j\}\in\mathcal{E}}c_{\{i,j\}}=N$. Similarly, by letting each node $i\in\mathcal{V}$ evenly distribute its observation $y_i$ to the $|\mathcal{N}_i|$ terms in $\sum_{\{i,j\}\in\mathcal{E}}c_{\{i,j\}}x_{\{i,j\}}(0)$, i.e.,
\begin{align}
x_{\{i,j\}}(0)=\frac{\frac{y_i}{|\mathcal{N}_i|}+\frac{y_j}{|\mathcal{N}_j|}}{c_{\{i,j\}}},\quad\forall\{i,j\}\in\mathcal{E},\label{eq:x=y|N|y|N|c}
\end{align}
we get $\sum_{\{i,j\}\in\mathcal{E}}c_{\{i,j\}}x_{\{i,j\}}(0)=\sum_{i\in\mathcal{V}}y_i$. Thus, \eqref{eq:c=1|N|1|N|} and \eqref{eq:x=y|N|y|N|c} together ensure \eqref{eq:sumcxsumc=sumyN}, which is necessary for achieving \eqref{eq:limx=x*}.

\begin{remark}\label{rem:parainitstat}
Obviously, \eqref{eq:c=1|N|1|N|} and \eqref{eq:x=y|N|y|N|c} are not the only way to select the $c_{\{i,j\}}$'s and $x_{\{i,j\}}(0)$'s. In fact, their selection may be posed as an optimization problem, analogous to the synchronous algorithms in \cite{XiaoL04, XiaoL07}. Nevertheless, \eqref{eq:c=1|N|1|N|} and \eqref{eq:x=y|N|y|N|c} have the virtue of being simple and inexpensive to implement: for every link $\{i,j\}\in\mathcal{E}$, both $c_{\{i,j\}}$ and $x_{\{i,j\}}(0)$ depend only on local information $|\mathcal{N}_i|$, $|\mathcal{N}_j|$, $y_i$, and $y_j$ that nodes $i$ and $j$ know, as opposed to on global information derived from the graph $\mathcal{G}$, which is typically difficult and costly to gather, but often the outcome of optimization.\hfill$\blacksquare$
\end{remark}

The system \eqref{eq:x=sumcxsumcuijx}, \eqref{eq:xh=sumcxsumc}, \eqref{eq:Pu=i=1N} with parameters \eqref{eq:c=1|N|1|N|} and initial states \eqref{eq:x=y|N|y|N|c} can be realized over the wireless network by having the nodes take the following actions: for every link $\{i,j\}\in\mathcal{E}$, nodes $i$ and $j$ each maintains a local copy of $x_{\{i,j\}}(k)$, denoted as $x_{ij}(k)$ and $x_{ji}(k)$, respectively, where they are meant to be always equal, so that the order of the subscripts is only used to indicate where they physically reside. Each node $i\in\mathcal{V}$, in addition to $x_{ij}(k)$ $\forall j\in\mathcal{N}_i$, also maintains $c_{\{i,j\}}$ $\forall j\in\mathcal{N}_i$ and $\hat{x}_i(k)$. To initialize the system, every node $i\in\mathcal{V}$ transmits $|\mathcal{N}_i|$ and $y_i$ each once, to every one-hop neighbor $j\in\mathcal{N}_i$, so that upon completion, each node $i\in\mathcal{V}$ can calculate $c_{\{i,j\}}$ $\forall j\in\mathcal{N}_i$ from \eqref{eq:c=1|N|1|N|}, $x_{ij}(0)$ $\forall j\in\mathcal{N}_i$ from \eqref{eq:x=y|N|y|N|c}, and $\hat{x}_i(0)$ from \eqref{eq:xh=sumcxsumc}. To evolve the system, at each iteration $k\in\mathbb{P}$, a node $u(k)\in\mathcal{V}$ is selected randomly and equiprobably based on \eqref{eq:Pu=i=1N} to initiate the iteration. To describe the subsequent actions, note that \eqref{eq:x=sumcxsumcuijx} and \eqref{eq:xh=sumcxsumc} imply: (i) $\hat{x}_{u(k)}(k)=\hat{x}_{u(k)}(k-1)$; (ii) $x_{u(k)j}(k)=\hat{x}_{u(k)}(k)$ $\forall j\in\mathcal{N}_{u(k)}$; (iii) $x_{ju(k)}(k)=\hat{x}_{u(k)}(k)$ $\forall j\in\mathcal{N}_{u(k)}$; (iv) $x_{j\ell}(k)=x_{j\ell}(k-1)$ $\forall\ell\in\mathcal{N}_j-\{u(k)\}$ $\forall j\in\mathcal{N}_{u(k)}$; (v) $\hat{x}_j(k)=\frac{\sum_{\ell\in\mathcal{N}_j}c_{\{j,\ell\}}x_{j\ell}(k)}{\sum_{\ell\in\mathcal{N}_j}c_{\{j,\ell\}}}$ $\forall j\in\mathcal{N}_{u(k)}$; (vi) $x_{\ell m}(k)=x_{\ell m}(k-1)$ $\forall m\in\mathcal{N}_\ell$ $\forall\ell\in\mathcal{V}-(\{u(k)\}\cup\mathcal{N}_{u(k)})$; and (vii) $\hat{x}_\ell(k)=\hat{x}_\ell(k-1)$ $\forall\ell\in\mathcal{V}-(\{u(k)\}\cup\mathcal{N}_{u(k)})$. To execute~(i) and~(ii), node $u(k)$, upon being selected to initiate iteration $k$, sets $\hat{x}_{u(k)}(k)$ and $x_{u(k)j}(k)$ $\forall j\in\mathcal{N}_{u(k)}$ all to $\hat{x}_{u(k)}(k-1)$. To execute~(iii), node $u(k)$ then transmits $\hat{x}_{u(k)}(k)$ once, to every one-hop neighbor $j\in\mathcal{N}_{u(k)}$, so that upon reception, each of them can set $x_{ju(k)}(k)$ to $\hat{x}_{u(k)}(k)$. Equations~(iv) and~(v) say that every neighbor $j\in\mathcal{N}_{u(k)}$ experiences no change in the rest of its local copies and, hence, can compute $\hat{x}_j(k)$ from~(v) upon finishing~(iii). Finally, (vi) and~(vii) say that the rest of the $N$ nodes, i.e., excluding node $u(k)$ and its one-hop neighbors, experience no change in the variables they maintain.

The above node actions define a distributed averaging algorithm that runs iteratively and asynchronously on the wireless network. We refer to this algorithm as {\em Random Hopwise Averaging} (RHA), since every iteration is {\em randomly} initiated and involves state variables associated with links within one {\em hop} of each other. RHA may be expressed in a compact algorithmic form as follows:

\begin{algorithm}[Random Hopwise Averaging]\label{alg:RHA}
\begin{algorithminit}{}
\item Each node $i\in\mathcal{V}$ transmits $|\mathcal{N}_i|$ and $y_i$ to every node $j\in\mathcal{N}_i$.
\item Each node $i\in\mathcal{V}$ creates variables $x_{ij}\in\mathbb{R}$ $\forall j\in\mathcal{N}_i$ and $\hat{x}_i\in\mathbb{R}$ and initializes them sequentially:
\\$x_{ij}\leftarrow\frac{\frac{y_i}{|\mathcal{N}_i|}+\frac{y_j}{|\mathcal{N}_j|}}{c_{\{i,j\}}},\quad\forall j\in\mathcal{N}_i,$
\\$\hat{x}_i\leftarrow\frac{\sum_{j\in\mathcal{N}_i}c_{\{i,j\}}x_{ij}}{\sum_{j\in\mathcal{N}_i}c_{\{i,j\}}}.$
\end{algorithminit}
\begin{algorithmoper}{At each iteration:}
\item A node, say, node $i$, is selected randomly and equiprobably out of the set $\mathcal{V}$ of $N$ nodes.
\item Node $i$ updates $x_{ij}$ $\forall j\in\mathcal{N}_i$:
\\$x_{ij}\leftarrow\hat{x}_i,\quad\forall j\in\mathcal{N}_i.$
\item Node $i$ transmits $\hat{x}_i$ to every node $j\in\mathcal{N}_i$.
\item Each node $j\in\mathcal{N}_i$ updates $x_{ji}$ and $\hat{x}_j$ sequentially:
\\$x_{ji}\leftarrow\hat{x}_i,$
\\$\hat{x}_j\leftarrow\frac{\sum_{\ell\in\mathcal{N}_j}c_{\{j,\ell\}}x_{j\ell}}{\sum_{\ell\in\mathcal{N}_j}c_{\{j,\ell\}}}.$
\end{algorithmoper}
\end{algorithm}

Observe from Algorithm~\ref{alg:RHA} that RHA requires an initialization overhead of $2N$ real-number transmissions to perform Step~1 (the $|\mathcal{N}_i|$'s are counted as real numbers, for simplicity). However, each iteration of RHA requires only transmission of a {\em single} message, consisting of exactly {\em one} real number, by the initiating node, in Step~5. Also notice that RHA fully exploits the broadcast nature of wireless medium, allowing everyone that hears the message to use it for revising their local variables, in Step~6. Therefore, RHA avoids issues~\ref{enu:overiter} and~\ref{enu:wastrece} with overlapping iterations and wasted receptions. Furthermore, as RHA operates asynchronously and calculates the average directly, it circumvents issues~\ref{enu:costdisc}--\ref{enu:compintequan} with costly discretization, clock synchronization, forced transmissions, and computing intermediate quantities. To show that it overcomes issues~\ref{enu:steastaterro} and~\ref{enu:lackconvguar} with steady-state errors and convergence guarantees, consider a quadratic Lyapunov function candidate $V:\mathbb{R}^L\rightarrow\mathbb{R}$, defined as
\begin{align}
V(\mathbf{x}(k))=\sum_{\{i,j\}\in\mathcal{E}}c_{\{i,j\}}(x_{\{i,j\}}(k)-x^*)^2.\label{eq:V=sumcxx*2}
\end{align}
Clearly, $V$ in \eqref{eq:V=sumcxx*2} is positive definite with respect to $(x^*,x^*,\ldots,x^*)\in\mathbb{R}^L$, and the condition
\begin{align}
\lim_{k\rightarrow\infty}V(\mathbf{x}(k))=0\label{eq:limV=0}
\end{align}
implies \eqref{eq:limx=x*} and thus \eqref{eq:limxh=x*}. The following lemma shows that $V(\mathbf{x}(k))$ is always non-increasing and quantifies its changes:

\begin{lemma}\label{lem:RHAVnonincr}
Consider the wireless network modeled in Section~\ref{sec:probform} and the use of RHA described in Algorithm~\ref{alg:RHA}. Then, for any sequence $(u(k))_{k=1}^\infty$, the sequence $(V(\mathbf{x}(k)))_{k=0}^\infty$ is non-increasing and satisfies
\begin{align}
V(\mathbf{x}(k))-V(\mathbf{x}(k-1))=-\sum_{j\in\mathcal{N}_{u(k)}}c_{\{u(k),j\}}(x_{\{u(k),j\}}(k-1)-\hat{x}_{u(k)}(k-1))^2,\quad\forall k\in\mathbb{P}.\label{eq:VV=sumcxxh2}
\end{align}
\end{lemma}

\begin{proof}
From \eqref{eq:V=sumcxx*2} and the bottom of \eqref{eq:x=sumcxsumcuijx}, $V(\mathbf{x}(k))-V(\mathbf{x}(k-1))=-\sum_{j\in\mathcal{N}_{u(k)}}c_{\{u(k),j\}}(-x_{\{u(k),j\}}^2(k)+2x_{\{u(k),j\}}(k)x^*+x_{\{u(k),j\}}^2(k-1)-2x_{\{u(k),j\}}(k-1)x^*)$ $\forall k\in\mathbb{P}$. Due to the top of \eqref{eq:x=sumcxsumcuijx}, the second term $-\sum_{j\in\mathcal{N}_{u(k)}}2c_{\{u(k),j\}}x_{\{u(k),j\}}(k)x^*$ cancels the fourth term $\sum_{j\in\mathcal{N}_{u(k)}}2c_{\{u(k),j\}}x_{\{u(k),j\}}(k-1)x^*$. Moreover, note from \eqref{eq:x=sumcxsumcuijx} and \eqref{eq:xh=sumcxsumc} that $x_{\{u(k),j\}}(k)=\hat{x}_{u(k)}(k-1)$ $\forall j\in\mathcal{N}_{u(k)}$. Hence, $V(\mathbf{x}(k))-V(\mathbf{x}(k-1))=-\sum_{j\in\mathcal{N}_{u(k)}}c_{\{u(k),j\}}(\hat{x}_{u(k)}^2(k-1)-2\hat{x}_{u(k)}(k-1)x_{\{u(k),j\}}(k)+x_{\{u(k),j\}}^2(k-1))$ $\forall k\in\mathbb{P}$. Due again to the top of \eqref{eq:x=sumcxsumcuijx}, the second term $\sum_{j\in\mathcal{N}_{u(k)}}2c_{\{u(k),j\}}\hat{x}_{u(k)}(k-1)x_{\{u(k),j\}}(k)$ equals $\sum_{j\in\mathcal{N}_{u(k)}}2c_{\{u(k),j\}}\hat{x}_{u(k)}(k-1)x_{\{u(k),j\}}(k-1)$. Thus, \eqref{eq:VV=sumcxxh2} holds. Since the right-hand side of \eqref{eq:VV=sumcxxh2} is nonpositive, $(V(\mathbf{x}(k)))_{k=0}^\infty$ is non-increasing.
\end{proof}

Lemma~\ref{lem:RHAVnonincr} says that $V(\mathbf{x}(k))\le V(\mathbf{x}(k-1))$ $\forall k\in\mathbb{P}$. Since $V(\mathbf{x}(k))\ge0$ $\forall\mathbf{x}(k)\in\mathbb{R}^L$, this implies that $\lim_{k\rightarrow\infty}V(\mathbf{x}(k))$ exists and is nonnegative. The following theorem asserts that this limit is almost surely zero, so that RHA is almost surely asymptotically convergent to $x^*$:

\begin{theorem}\label{thm:RHAasymconv}
Consider the wireless network modeled in Section~\ref{sec:probform} and the use of RHA described in Algorithm~\ref{alg:RHA}. Then, with probability $1$, \eqref{eq:limV=0}, \eqref{eq:limx=x*}, and \eqref{eq:limxh=x*} hold.
\end{theorem}

\begin{proof}
By associating the line graph of $\mathcal{G}$ with the graph in \cite{Fagnani08}, RHA may be viewed as a special case of the algorithm~(1) in \cite{Fagnani08}. Note from \eqref{eq:x=sumcxsumcuijx} and \eqref{eq:c=1|N|1|N|} that the diagonal entries of $\mathbf{A}_i$ $\forall i\in\mathcal{V}$ are positive, from \eqref{eq:Pu=i=1N} that $\operatorname{P}\{\mathbf{A}_{u(k)}=\mathbf{A}_i\}=\frac{1}{N}$ $\forall k\in\mathbb{P}$ $\forall i\in\mathcal{V}$, and from the connectedness of $\mathcal{G}$ that its line graph is connected. Thus, by Corollary~3.2 of \cite{Fagnani08}, with probability $1$, $\exists\tilde{x}\in\mathbb{R}$ such that $\lim_{k\rightarrow\infty}x_{\{i,j\}}(k)=\tilde{x}$ $\forall\{i,j\}\in\mathcal{E}$. Due to \eqref{eq:x*=1Nsumy}, \eqref{eq:sumcx=sumcx}, and \eqref{eq:sumcxsumc=sumyN}, $\tilde{x}=x^*$, i.e., \eqref{eq:limx=x*} holds almost surely. Because of \eqref{eq:V=sumcxx*2} and \eqref{eq:xh=sumcxsumc}, so do \eqref{eq:limV=0} and \eqref{eq:limxh=x*}.
\end{proof}

As it follows from Theorem~\ref{thm:RHAasymconv} and the above, RHA solves the distributed averaging problem, while eliminating deficiencies~\ref{enu:costdisc}--\ref{enu:lackconvguar} facing the existing schemes except for~\ref{enu:uncoiter}, on lack of control. Lemma~\ref{lem:RHAVnonincr} above also says that $V$ in \eqref{eq:V=sumcxx*2} is a {\em common} quadratic Lyapunov function for the linear switched system \eqref{eq:x=Aux}. This $V$ will be used next to introduce control and remove~\ref{enu:uncoiter}.

\section{Controlled Hopwise Averaging}\label{sec:CHA}

\subsection{Motivation for Feedback Iteration Control}\label{ssec:motifeeditercont}

RHA operates by executing \eqref{eq:x=sumcxsumcuijx} or \eqref{eq:x=Aux} according to $(u(k))_{k=1}^\infty$. Although, by Theorem~\ref{thm:RHAasymconv}, almost any $(u(k))_{k=1}^\infty$ can drive all the $\hat{x}_i(k)$'s in \eqref{eq:xh=sumcxsumc} to any neighborhood of $x^*$, certain sequences require fewer iterations (and, hence, fewer real-number transmissions) to do so than others, yielding better bandwidth/energy efficiency. To see this, consider the following proposition:

\begin{proposition}\label{pro:RHAidemcomm}
The matrices $\mathbf{A}_1,\mathbf{A}_2,\ldots,\mathbf{A}_N$ in \eqref{eq:x=Aux} are idempotent, i.e., $\mathbf{A}_i^2=\mathbf{A}_i$ $\forall i\in\mathcal{V}$. Moreover, $\mathbf{A}_i$ and $\mathbf{A}_j$ are commutative whenever $\{i,j\}\notin\mathcal{E}$, i.e., $\mathbf{A}_i\mathbf{A}_j=\mathbf{A}_j\mathbf{A}_i$ $\forall i,j\in\mathcal{V}$, $\{i,j\}\notin\mathcal{E}$.
\end{proposition}

\begin{proof}
Notice from \eqref{eq:x=sumcxsumcuijx} and \eqref{eq:x=Aux} that for any $i\in\mathcal{V}$, if $\mathbf{x}(k)=\mathbf{A}_i\mathbf{x}(k-1)$, then $x_{\{i,j\}}(k)$ $\forall j\in\mathcal{N}_i$ are set equal to the same convex combination of $x_{\{i,j\}}(k-1)$ $\forall j\in\mathcal{N}_i$, and $x_{\{p,q\}}(k)=x_{\{p,q\}}(k-1)$ $\forall\{p,q\}\in\mathcal{E}-\cup_{j\in\mathcal{N}_i}\{\{i,j\}\}$. Thus, $\mathbf{A}_i\mathbf{x}(k)=\mathbf{x}(k)$, so that $\mathbf{A}_i^2=\mathbf{A}_i$. Moreover, for any $i,j\in\mathcal{V}$ with $\{i,j\}\notin\mathcal{E}$, because $\{\{i,\ell\}:\ell\in\mathcal{N}_i\}\cap\{\{j,\ell\}:\ell\in\mathcal{N}_j\}=\emptyset$, $\mathbf{A}_i\mathbf{A}_j=\mathbf{A}_j\mathbf{A}_i$.
\end{proof}

The idempotence and partial commutativity of $\mathbf{A}_1,\mathbf{A}_2,\ldots,\mathbf{A}_N$ from Proposition~\ref{pro:RHAidemcomm}, together with the fact that the switched system \eqref{eq:x=Aux} may be stated as $\mathbf{x}(k)=\mathbf{A}_{u(k)}\mathbf{A}_{u(k-1)}\cdots\mathbf{A}_{u(1)}\mathbf{x}(0)$ $\forall k\in\mathbb{P}$, imply that for a given $(u(k))_{k=1}^\infty$, the event $\mathbf{x}(k)=\mathbf{x}(k-1)$ can occur for quite a few $k$'s, each of which signifies a wasted iteration. Furthermore, if the event $\mathbf{x}(k)=\mathbf{x}(k-1)$ does occur for at least one $k$, then by deleting from $(u(k))_{k=1}^\infty$ some of its elements that correspond to the wasted iterations, we obtain a new sequence $(u'(k))_{k=1}^\infty$ that is more efficient. To illustrate these two points, consider, for instance, a $5$-node cycle graph with $\mathcal{V}=\{1,2,3,4,5\}$ and $\mathcal{E}=\{\{1,2\},\{2,3\},\{3,4\},\{4,5\},\{5,1\}\}$. Notice that if $(u(k))_{k=1}^\infty=(1,\underline{1},3,4,\underline{1},2,\underline{4},5,\underline{2},\underline{5},\ldots)$, then as many as $5$ out of the first $10$ iterations---namely, those underlined elements---are wasted. By deleting these underlined elements and keeping the rest intact, we obtain a new sequence $(u'(k))_{k=1}^\infty=(1,3,4,2,5,\ldots)$ that is $5$ real-number transmissions more efficient than $(u(k))_{k=1}^\infty$.

The preceding analysis shows that RHA is prone to wasteful iterations, which is a primary reason why certain sequences are more efficient than others. RHA, however, makes no attempt to distinguish the sequences, as it lets every possible $(u(k))_{k=1}^\infty$ be equiprobable, via \eqref{eq:Pu=i=1N}. In other words, it does not try to {\em control} how the asynchronous iterations occur and, thus, suffers from~\ref{enu:uncoiter}.

\begin{remark}\label{rem:RHAidemcomm}
Wasteful iterations incurred by idempotent and partially commutative operations are not an attribute unique to RHA, but one that is shared by Pairwise Averaging \cite{Tsitsiklis84}, Anti-Entropy Aggregation \cite{Jelasity04, Montresor04}, Randomized Gossip Algorithm \cite{Boyd06}, and Distributed Random Grouping \cite{ChenJY06} (indeed, the examples provided in~\ref{enu:uncoiter} against the latter two algorithms were created from this attribute). What is different is that in this paper, we view the attribute as a limitation and find ways to overcome it, whereas in \cite{Tsitsiklis84, Jelasity04, Montresor04, Boyd06, ChenJY06}, the attribute was not viewed as such.\hfill$\blacksquare$
\end{remark}

One way to control the iterations, alluded to in Remark~\ref{rem:alterandequi}, is to replace \eqref{eq:Pu=i=1N} with a general distribution $\operatorname{P}\{u(k)=i\}=p_i$ $\forall k\in\mathbb{P}$ $\forall i\in\mathcal{V}$ and then choose the $p_i$'s to maximize efficiency, before any averaging task begins. This approach, however, has an inherent shortcoming: because the $p_i$'s are optimized once-and-for-all, they are constant and do not adapt to $\mathbf{x}(k)$ during runtime. Hence, optimal or not, the $p_i$'s almost surely would produce inefficient, wasteful $(u(k))_{k=1}^\infty$. The fact that the nodes do not adjust the $p_i$'s based on information they pick up during runtime also suggests that this way of controlling the iterations may be considered {\em open loop}.

The aforementioned shortcoming of open-loop iteration control raises the question of whether it is possible to introduce some form of {\em closed-loop} iteration control as a means to generate efficient, non-wasteful $(u(k))_{k=1}^\infty$. Obviously, to carry out closed-loop iteration control, feedback is needed. Due to the distributed nature of the network, however, feedback may be expensive to acquire: if an algorithm demands that the feedback used by a node be a function of state variables maintained by other nodes, then additional communications are necessary to implement the feedback. Such communications can produce plenty of real-number transmissions, which must all count toward the total real-number transmissions, when evaluating the algorithm's bandwidth/energy efficiency. Thus, in the design of feedback algorithms, the cost of ``closing the loop'' cannot be overlooked.

In this section, we first describe an approach to closed-loop iteration control, which leads to highly efficient and surely non-wasteful $(u(k))_{k=1}^\infty$ at {\em zero} feedback cost. Based on this approach, we then present and analyze two modified versions of RHA: an ideal version and a practical one.

\subsection{Approach to Feedback Iteration Control}\label{ssec:apprfeeditercont}

Note that with RHA, $(u(k))_{k=1}^\infty$ is undefined at the moment an averaging task begins and is gradually defined, one element per iteration, as time elapses, i.e., when a node $i\in\mathcal{V}$ initiates an iteration $k\in\mathbb{P}$, the element $u(k)$ becomes defined and is given by $u(k)=i$. Thus, by controlling {\em when} to initiate an iteration, the nodes may jointly shape the value of $(u(k))_{k=1}^\infty$. With RHA, this opportunity to shape $(u(k))_{k=1}^\infty$ is not utilized, as the nodes simply randomly and equiprobably decide when to initiate an iteration. To exploit the opportunity, suppose henceforth that the nodes wish to control when to initiate an iteration using some form of {\em feedback}. The questions are:
\begin{enumerate}
\renewcommand{\theenumi}{Q\arabic{enumi}}\itemsep-\parsep
\item What feedback to use, so that the corresponding feedback cost is minimal?\label{enu:whatfeedcostmini}
\item How to control, so that the resulting $(u(k))_{k=1}^\infty$ is highly efficient?\label{enu:howconthigheffi}
\item How to control, so that the resulting $(u(k))_{k=1}^\infty$ is surely non-wasteful?\label{enu:howcontsurenonwast}
\end{enumerate}

To answer questions~\ref{enu:whatfeedcostmini}--\ref{enu:howcontsurenonwast}, we first show that RHA, along with the common quadratic Lyapunov function $V$ of \eqref{eq:V=sumcxx*2}, exhibits the following features:
\begin{enumerate}
\renewcommand{\theenumi}{F\arabic{enumi}}\itemsep-\parsep
\item Although the nodes never know the value of $V$, every one of them at any time knows by how much the value would drop if it suddenly initiates an iteration.\label{enu:knowVdrop}
\item The faster $(u(k))_{k=1}^\infty$ makes the value of $V$ drop to zero, the more efficient it is.\label{enu:Vdropfastueffi}
\item If the value of $V$ does not drop after an iteration, then the iteration is wasted, causing $(u(k))_{k=1}^\infty$ to be wasteful. The converse is also true.\label{enu:Vdrop0uwast}
\end{enumerate}

The first part of feature~\ref{enu:knowVdrop} can be seen by noting that $V(\mathbf{x}(k))$ in \eqref{eq:V=sumcxx*2} depends on $c_{\{i,j\}}$ $\forall\{i,j\}\in\mathcal{E}$, $x_{\{i,j\}}(k)$ $\forall\{i,j\}\in\mathcal{E}$, and $x^*$, whereas each node $i\in\mathcal{V}$ only knows $c_{\{i,j\}}$ $\forall j\in\mathcal{N}_i$ and $x_{\{i,j\}}(k)$ $\forall j\in\mathcal{N}_i$. To see the second part, suppose a node $i\in\mathcal{V}$ initiates an iteration $k\in\mathbb{P}$ at some time instant $t$, so that $u(k)=i$ by definition. Observe from Lemma~\ref{lem:RHAVnonincr} that whoever node $u(k)$ is, upon completing this iteration, the value of $V$ would drop from $V(\mathbf{x}(k-1))$ to $V(\mathbf{x}(k))$ by an amount equal to the right-hand side of \eqref{eq:VV=sumcxxh2}. To compactly represent this drop, for each $i\in\mathcal{V}$ let $\Delta V_i:\mathbb{R}^L\rightarrow\mathbb{R}$ be a positive semidefinite quadratic function, defined as
\begin{align}
\Delta V_i(\mathbf{x}(k))=\sum_{j\in\mathcal{N}_i}c_{\{i,j\}}(x_{\{i,j\}}(k)-\hat{x}_i(k))^2,\quad\forall k\in\mathbb{N},\label{eq:DV=sumcxxh2}
\end{align}
where $\hat{x}_i(k)$ is as in \eqref{eq:xh=sumcxsumc}. Then, with \eqref{eq:DV=sumcxxh2}, \eqref{eq:VV=sumcxxh2} may be written as
\begin{align}
V(\mathbf{x}(k))-V(\mathbf{x}(k-1))=-\Delta V_{u(k)}(\mathbf{x}(k-1)),\quad\forall k\in\mathbb{P},\label{eq:VV=DVu}
\end{align}
where $\Delta V_{u(k)}(\mathbf{x}(k-1))$ in \eqref{eq:VV=DVu} represents the amount of drop, i.e.,
\begin{align}
\Delta V_{u(k)}(\mathbf{x}(k-1))=\sum_{j\in\mathcal{N}_{u(k)}}c_{\{u(k),j\}}(x_{\{u(k),j\}}(k-1)-\hat{x}_{u(k)}(k-1))^2,\quad\forall k\in\mathbb{P}.\label{eq:DVu=sumcxxh2}
\end{align}
Notice that $\Delta V_{u(k)}(\mathbf{x}(k-1))$ in \eqref{eq:DVu=sumcxxh2} depends on parameters and variables maintained by node $u(k)$, whose values are known to node $u(k)$ prior to iteration $k$ at time $t$. Therefore, before initiating this iteration at time $t$, node $u(k)$ already knows that the value of $V$ would drop by $\Delta V_{u(k)}(\mathbf{x}(k-1))$. Since $t$, $k$, and $u(k)$ are arbitrary, this means that every node $i\in\mathcal{V}$ at any time knows by how much the value of $V$ would drop if it suddenly initiates an iteration (i.e., by $\Delta V_i(\mathbf{x}(\cdot))$). This establishes feature~\ref{enu:knowVdrop}. To show feature~\ref{enu:Vdropfastueffi}, recall that: (i) $V(\mathbf{x}(k))$ in \eqref{eq:V=sumcxx*2} is a measure of the deviation of the $x_{\{i,j\}}(k)$'s from $x^*$; (ii) the $\hat{x}_i(k)$'s in \eqref{eq:xh=sumcxsumc} are convex combinations of the $x_{\{i,j\}}(k)$'s; (iii) bandwidth/energy efficiency is measured by the number of real-number transmissions needed for all the $\hat{x}_i(k)$'s to converge to a given neighborhood of $x^*$; and (iv) RHA in Algorithm~\ref{alg:RHA} has a fixed, one real-number transmission per iteration. Hence, the faster $(u(k))_{k=1}^\infty$ drives $V(\mathbf{x}(k))$ to zero, the faster it drives the $x_{\{i,j\}}(k)$'s and $\hat{x}_i(k)$'s to $x^*$ (due to~(i) and~(ii)), and the more efficient it is (due to~(iii) and~(iv)). Finally, to show feature~\ref{enu:Vdrop0uwast}, suppose $V(\mathbf{x}(k))=V(\mathbf{x}(k-1))$ after an iteration $k\in\mathbb{P}$. Then, it follows from \eqref{eq:VV=DVu} that $\Delta V_{u(k)}(\mathbf{x}(k-1))=0$, from \eqref{eq:DVu=sumcxxh2} that $x_{\{u(k),j\}}(k-1)$ $\forall j\in\mathcal{N}_{u(k)}$ are equal, and from \eqref{eq:x=sumcxsumcuijx} that $\mathbf{x}(k)=\mathbf{x}(k-1)$. Thus, iteration $k$ is wasted. The converse is also true, as $\mathbf{x}(k)=\mathbf{x}(k-1)$ implies $V(\mathbf{x}(k))=V(\mathbf{x}(k-1))$.

Having demonstrated features~\ref{enu:knowVdrop}--\ref{enu:Vdrop0uwast}, we now use them to answer questions~\ref{enu:whatfeedcostmini}--\ref{enu:howcontsurenonwast}. Feature~\ref{enu:knowVdrop} suggests that every node $i\in\mathcal{V}$ may use $\Delta V_i(\mathbf{x}(\cdot))$, which it always knows, as feedback to control, on its own, when to initiate an iteration. As the feedbacks $\Delta V_i(\mathbf{x}(\cdot))$'s are locally available and the control decisions are made locally, the resulting feedback control architecture is fully decentralized, requiring zero communication cost to realize. Therefore, an answer to question~\ref{enu:whatfeedcostmini} is:
\begin{enumerate}
\renewcommand{\theenumi}{A\arabic{enumi}}\itemsep-\parsep
\item Each node $i\in\mathcal{V}$ uses $\Delta V_i(\mathbf{x}(\cdot))$ as feedback to control when to initiate an iteration.\label{enu:useDVfeedcontwheninit}
\end{enumerate}
Feature~\ref{enu:Vdropfastueffi} suggests that, to produce highly efficient $(u(k))_{k=1}^\infty$, the nodes may focus on making the value of $V$ drop significantly after each iteration, especially initially. In other words, they may focus on letting every iteration be initiated by a node $i$ with a relatively large $\Delta V_i(\mathbf{x}(\cdot))$. With architecture~\ref{enu:useDVfeedcontwheninit}, this may be accomplished if nodes with larger $\Delta V_i(\mathbf{x}(\cdot))$'s would rush to initiate, while nodes with smaller $\Delta V_i(\mathbf{x}(\cdot))$'s would wait longer. Hence, an answer to question~\ref{enu:howconthigheffi} is:
\begin{enumerate}
\renewcommand{\theenumi}{A\arabic{enumi}}\itemsep-\parsep\setcounter{enumi}{1}
\item The larger $\Delta V_i(\mathbf{x}(\cdot))$ is, the sooner node $i$ initiates an iteration (i.e., the smaller $\Delta V_i(\mathbf{x}(\cdot))$ is, the longer node $i$ waits).\label{enu:largDVsooninit}
\end{enumerate}
Finally, feature~\ref{enu:Vdrop0uwast} suggests that, to generate surely non-wasteful $(u(k))_{k=1}^\infty$, the value of $V$ must strictly decrease after each iteration. With architecture~\ref{enu:useDVfeedcontwheninit}, this can be achieved if nodes with zero $\Delta V_i(\mathbf{x}(\cdot))$'s would refrain from initiating an iteration. Thus, an answer to question~\ref{enu:howcontsurenonwast} is:
\begin{enumerate}
\renewcommand{\theenumi}{A\arabic{enumi}}\itemsep-\parsep\setcounter{enumi}{2}
\item Whenever $\Delta V_i(\mathbf{x}(\cdot))=0$, node $i$ refrains from initiating an iteration.\label{enu:whenDV0refrinit}
\end{enumerate}

Answers~\ref{enu:useDVfeedcontwheninit}--\ref{enu:whenDV0refrinit} describe a greedy, decentralized approach to feedback iteration control, where potential drops $\Delta V_i(\mathbf{x}(\cdot))$'s in the value of $V$ are used to drive the asynchronous iterations. This approach may be viewed as a greedy approach because the nodes seek to make the value of $V$ drop as much as possible at each iteration, without considering the future. Because the nodes also seek to fully exploit the broadcast nature of every wireless transmission (a feature inherited from Steps~5 and~6 of RHA), this approach strives to ``make the most'' out of each iteration. Note that although Lyapunov functions have been used to analyze distributed averaging and consensus algorithms (e.g., in the form of a disagreement function \cite{Olfati-Saber04} or a set-valued convex hull \cite{Moreau05}), their use for {\em controlling} such algorithms has not been reported. Therefore, this approach represents a new way to apply Lyapunov stability theory.

\subsection{Ideal Version}\label{ssec:ideavers}

In this subsection, we use the aforementioned approach to create an ideal, modified version of RHA, which possesses strong convergence properties that motivate a practical version.

The above approach wants the nodes to try to be greedy. Thus, it is of interest to analyze an ideal scenario where, instead of just trying, the nodes actually succeed at being greedy, ensuring that every iteration $k\in\mathbb{P}$ is initiated by a node $i\in\mathcal{V}$ with the maximum $\Delta V_i(\mathbf{x}(k-1))$, i.e.,
\begin{align}
u(k)\in\operatornamewithlimits{arg\,max}_{i\in\mathcal{V}}\Delta V_i(\mathbf{x}(k-1)),\quad\forall k\in\mathbb{P},\label{eq:u=argmaxDV}
\end{align}
so that $V(\mathbf{x}(k-1))$ drops maximally to $V(\mathbf{x}(k))$ for every $k\in\mathbb{P}$. Notice that \eqref{eq:u=argmaxDV} does not always uniquely determine $u(k)$: when multiple nodes have the same maximum, $u(k)$ may be any of these nodes. Although $u(k)$ can be made unique (e.g., by letting $u(k)$ be the minimum of $\operatornamewithlimits{arg\,max}_{i\in\mathcal{V}}\Delta V_i(\mathbf{x}(k-1))$), in the analysis below we will allow for arbitrary $u(k)$ satisfying \eqref{eq:u=argmaxDV}. Also note that in the rare case where $\Delta V_i(\mathbf{x}(k^*-1))=0$ $\forall i\in\mathcal{V}$ for some $k^*\in\mathbb{P}$, due to \eqref{eq:x*=1Nsumy}, \eqref{eq:sumcx=sumcx}, \eqref{eq:sumcxsumc=sumyN}, \eqref{eq:DV=sumcxxh2}, and the connectedness of the graph $\mathcal{G}$, we have $x_{\{i,j\}}(k^*-1)=x^*$ $\forall\{i,j\}\in\mathcal{E}$ and $\hat{x}_i(k^*-1)=x^*$ $\forall i\in\mathcal{V}$, thereby solving the problem in finite time. Furthermore, due to~\ref{enu:whenDV0refrinit}, all the nodes would refrain from initiating iteration $k^*$ (and beyond), thereby terminating the algorithm in finite time and causing $x_{\{i,j\}}(k)$ $\forall\{i,j\}\in\mathcal{E}$, $\hat{x}_i(k)$ $\forall i\in\mathcal{V}$, $u(k)$, and $V(\mathbf{x}(k))$ to be undefined $\forall k\ge k^*$. In the analysis below, however, we will allow the algorithm to keep executing according to \eqref{eq:u=argmaxDV}, so that $x_{\{i,j\}}(k)$ $\forall\{i,j\}\in\mathcal{E}$, $\hat{x}_i(k)$ $\forall i\in\mathcal{V}$, $u(k)$, and $V(\mathbf{x}(k))$ are defined $\forall k$.

Equation \eqref{eq:u=argmaxDV}, together with \eqref{eq:x=sumcxsumcuijx}, \eqref{eq:xh=sumcxsumc}, \eqref{eq:c=1|N|1|N|}, \eqref{eq:x=y|N|y|N|c}, and \eqref{eq:DV=sumcxxh2}, defines a networked dynamical system that switches among $N$ different dynamics, depending on where the state is in the state space, i.e., if $\mathbf{x}(k-1)$ is such that $\Delta V_i(\mathbf{x}(k-1))>\Delta V_j(\mathbf{x}(k-1))$ $\forall j\in\mathcal{V}-\{i\}$, then $\mathbf{x}(k)=\mathbf{A}_i\mathbf{x}(k-1)$. This system may be expressed in the form of an algorithm---which we refer to as {\em Ideal Controlled Hopwise Averaging} (ICHA)---as follows:

\begin{algorithm}[Ideal Controlled Hopwise Averaging]\label{alg:ICHA}
\begin{algorithminit}{}
\item Each node $i\in\mathcal{V}$ transmits $|\mathcal{N}_i|$ and $y_i$ to every node $j\in\mathcal{N}_i$.
\item Each node $i\in\mathcal{V}$ creates variables $x_{ij}\in\mathbb{R}$ $\forall j\in\mathcal{N}_i$, $\hat{x}_i\in\mathbb{R}$, and $\Delta V_i\in[0,\infty)$ and initializes them sequentially:
\\$x_{ij}\leftarrow\frac{\frac{y_i}{|\mathcal{N}_i|}+\frac{y_j}{|\mathcal{N}_j|}}{c_{\{i,j\}}},\quad\forall j\in\mathcal{N}_i,$
\\$\hat{x}_i\leftarrow\frac{\sum_{j\in\mathcal{N}_i}c_{\{i,j\}}x_{ij}}{\sum_{j\in\mathcal{N}_i}c_{\{i,j\}}},$
\\$\Delta V_i\leftarrow\sum_{j\in\mathcal{N}_i}c_{\{i,j\}}(x_{ij}-\hat{x}_i)^2.$
\end{algorithminit}
\begin{algorithmoper}{At each iteration:}
\item Let $i\in\operatornamewithlimits{arg\,max}_{j\in\mathcal{V}}\Delta V_j$.
\item Node $i$ updates $x_{ij}$ $\forall j\in\mathcal{N}_i$ and $\Delta V_i$ sequentially:
\\$x_{ij}\leftarrow\hat{x}_i,\quad\forall j\in\mathcal{N}_i,$
\\$\Delta V_i\leftarrow0.$
\item Node $i$ transmits $\hat{x}_i$ to every node $j\in\mathcal{N}_i$.
\item Each node $j\in\mathcal{N}_i$ updates $x_{ji}$, $\hat{x}_j$, and $\Delta V_j$ sequentially:
\\$x_{ji}\leftarrow\hat{x}_i,$
\\$\hat{x}_j\leftarrow\frac{\sum_{\ell\in\mathcal{N}_j}c_{\{j,\ell\}}x_{j\ell}}{\sum_{\ell\in\mathcal{N}_j}c_{\{j,\ell\}}},$
\\$\Delta V_j\leftarrow\sum_{\ell\in\mathcal{N}_j}c_{\{j,\ell\}}(x_{j\ell}-\hat{x}_j)^2.$
\end{algorithmoper}
\end{algorithm}

Algorithm~\ref{alg:ICHA}, or ICHA, is identical to RHA in Algorithm~\ref{alg:RHA} except that each node $i$ also maintains $\Delta V_i$, in Steps~2, 4, and~6, and that each iteration is initiated by a node $i$ experiencing the maximum $\Delta V_i$, in Step~3. Note that ``$\Delta V_i\leftarrow0$'' in Step~4 is equivalent to ``$\Delta V_i\leftarrow\sum_{j\in\mathcal{N}_i}c_{\{i,j\}}(x_{ij}-\hat{x}_i)^2$'' since $x_{ij}$ $\forall j\in\mathcal{N}_i$ and $\hat{x}_i$ are equal at that point. The fact that $\Delta V_i$ goes from being the maximum to zero whenever node $i$ initiates an iteration also suggests that it may be a while before $\Delta V_i$ becomes the maximum again, causing node $i$ to initiate another iteration.

The convergence properties of ICHA on general networks are characterized in the following theorem, in which $\mathbf{1}_n\in\mathbb{R}^n$ and $\hat{\mathbf{x}}(k)\in\mathbb{R}^N$ denote, respectively, the vectors obtained by stacking $n$ $1$'s and the $N$ $\hat{x}_i(k)$'s:

\begin{theorem}\label{thm:ICHAexpconvgene}
Consider the wireless network modeled in Section~\ref{sec:probform} and the use of ICHA described in Algorithm~\ref{alg:ICHA}. Then,
\begin{align}
V(\mathbf{x}(k))&\le(1-\tfrac{1}{\gamma})V(\mathbf{x}(k-1)),\quad\forall k\in\mathbb{P},\label{eq:V<=11gV}\displaybreak[0]\\
\|\mathbf{x}(k)-x^*\mathbf{1}_L\|&\le\sqrt{\tfrac{V(\mathbf{x}(0))\max_{i\in\mathcal{V}}|\mathcal{N}_i|}{2}}(1-\tfrac{1}{\gamma})^{k/2},\quad\forall k\in\mathbb{N},\label{eq:||xx*1||<=sqrtVmax|N|211gk2}\displaybreak[0]\\
\|\hat{\mathbf{x}}(k)-x^*\mathbf{1}_N\|&\le\sqrt{\tfrac{2V(\mathbf{x}(0))\max_{i\in\mathcal{V}}|\mathcal{N}_i|}{\min_{i\in\mathcal{V}}|\mathcal{N}_i|+\max_{i\in\mathcal{V}}|\mathcal{N}_i|}}(1-\tfrac{1}{\gamma})^{k/2},\quad\forall k\in\mathbb{N},\label{eq:||xhx*1||<=sqrt2Vmax|N|min|N|max|N|11gk2}
\end{align}
where $\gamma\in[\frac{N}{2}+1,N^3-2N^2+\frac{N}{2}+1]$ is given by
\begin{align}
\gamma=\frac{N}{2}+\alpha+\frac{(N^2-\beta)(3(N-1)-D)(D+1)}{2N},\label{eq:g=N2aN2b3N1DD12N}
\end{align}
and where $\alpha=\max_{\{i,j\}\in\mathcal{E}}\frac{b_i+b_j}{c_{\{i,j\}}}\in[1,\frac{N^2-2N+2}{2}]$, $\beta=\sum_{i\in\mathcal{V}}\sum_{j\in\mathcal{N}_i\cup\{i\}}b_ib_j\in[N+\frac{L}{2}(1+\frac{1}{N-1})^2,N^2]$, $b_i=\frac{1}{2}\sum_{j\in\mathcal{N}_i}c_{\{i,j\}}$ $\forall i\in\mathcal{V}$, and $D$ is the network diameter.
\end{theorem}

\begin{proof}
See Appendix~\ref{ssec:proofthmICHAexpconvgene}.
\end{proof}

Theorem~\ref{thm:ICHAexpconvgene} says that ICHA is exponentially convergent on any network, ensuring that $V(\mathbf{x}(k))$, $\|\mathbf{x}(k)-x^*\mathbf{1}_L\|$, and $\|\hat{\mathbf{x}}(k)-x^*\mathbf{1}_N\|$ all go to zero exponentially fast, at a rate that is no worse than $1-\frac{1}{\gamma}$ or $(1-\frac{1}{\gamma})^{1/2}$, so that $\gamma$ in \eqref{eq:g=N2aN2b3N1DD12N} represents a bound on the convergence rate. It also says that the bound $\gamma$ is between $\Omega(N)$ and $O(N^3)$ and depends only on $N$, $D$, and the $|\mathcal{N}_i|$'s, making it easy to compute. The following corollary lists the bound $\gamma$ for a number of common graphs:

\begin{corollary}\label{cor:ICHAexpconvgene}
The constant $\gamma$ in \eqref{eq:g=N2aN2b3N1DD12N} becomes:
\begin{enumerate}
\renewcommand{\theenumi}{G\arabic{enumi}}\itemsep-\parsep
\item $\gamma=N^3-4N^2+\frac{9}{2}N+\frac{5}{4}$ for a path graph with $N\ge5$,\label{enu:genepath}
\item $\gamma=\frac{5}{8}N^3-\frac{15}{8}N^2-\frac{1}{8}N+\frac{31}{8}$ if $N$ is odd and $\gamma=\frac{5}{8}N^3-\frac{11}{8}N^2-\frac{5}{2}N+\frac{13}{2}$ if $N$ is even for a cycle graph,\label{enu:genecycl}
\item $\gamma=\frac{N}{2}+K+\frac{(N-K-1)(3(N-1)-D)(D+1)}{2}$ for a $K$-regular graph with $K\ge2$,\label{enu:generegu}
\item $\gamma=\frac{3}{2}N-1$ for a complete graph.\label{enu:genecomp}
\end{enumerate}
\end{corollary}

\begin{proof}
For a path graph with $N\ge5$, $\alpha=\frac{9}{4}$, $\beta=3N-1$, and $D=N-1$. For a cycle graph, $\alpha=2$, $\beta=3N$, $D=\frac{N-1}{2}$ if $N$ is odd, and $D=\frac{N}{2}$ if $N$ is even. For a $K$-regular graph with $K\ge2$, $\alpha=K$ and $\beta=N(K+1)$. For a complete graph, $\alpha=N-1$ and $\beta=N^2$. Hence, \ref{enu:genepath}--\ref{enu:genecomp} hold.
\end{proof}

Each bound $\gamma$ in Corollary~\ref{cor:ICHAexpconvgene} is obtained by specializing \eqref{eq:g=N2aN2b3N1DD12N} for arbitrary graphs to a specific one. Conceivably, tighter bounds may be obtained by working with each of these graphs individually, exploiting their particular structure. Theorem~\ref{thm:ICHAexpconvspec} below shows that this is indeed the case with path and cycle graphs ($6$ and $15$ times tighter, respectively), besides providing additional bounds for regular and strongly regular graphs:

\begin{theorem}\label{thm:ICHAexpconvspec}
Consider the wireless network modeled in Section~\ref{sec:probform} and the use of ICHA described in Algorithm~\ref{alg:ICHA}. Then, \eqref{eq:V<=11gV}--\eqref{eq:||xhx*1||<=sqrt2Vmax|N|min|N|max|N|11gk2} hold with:
\begin{enumerate}
\renewcommand{\theenumi}{S\arabic{enumi}}\itemsep-\parsep
\item $\gamma=\frac{N^3}{6}-\frac{13}{6}N+3$ for a path graph with $N\ge4$,\label{enu:specpath}
\item $\gamma=\frac{N^3}{24}+\frac{7}{12}N-2+\frac{11}{8N}$ if $N$ is odd and $\gamma=\frac{N^3}{24}+\frac{5}{6}N-3+\frac{4}{N}$ if $N$ is even for a cycle graph,\label{enu:speccycl}
\item $\gamma=\frac{N}{2}+K+\frac{KD(D+1)(N-K-1)}{2}$ for a $K$-regular graph with $K\ge2$,\label{enu:specregu}
\item $\gamma=\frac{N}{2}+K+\frac{K(\mu+2)(N-K-1)}{\mu}$ for a $(N,K,\lambda,\mu)$-strongly regular graph with $\mu\ge1$.\label{enu:specstro}
\end{enumerate}
\end{theorem}

\begin{proof}
See Appendix~\ref{ssec:proofthmICHAexpconvspec}.
\end{proof}

\begin{figure}[tb]
\centering\includegraphics{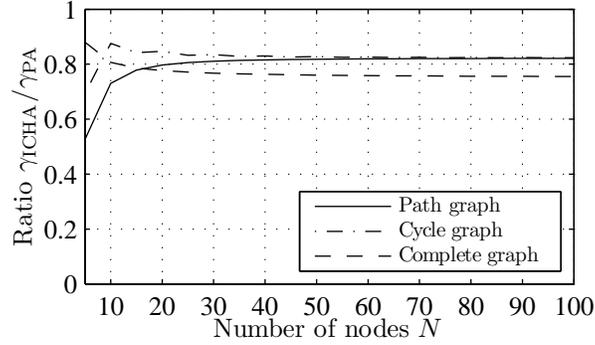}
\caption{Comparison between the stochastic convergence rate $1-\frac{1}{\gamma_{\text{PA}}}$ of PA and the deterministic bound $1-\frac{1}{\gamma_{\text{ICHA}}}$ on convergence rate of ICHA for path, cycle, and complete graphs.}
\label{fig:gammagamma_N}
\end{figure}

Recently, \cite{Fagnani08} studied, among other things, the convergence rate of Pairwise Averaging (PA) \cite{Tsitsiklis84}. The results in \cite{Fagnani08} are different from those above in three notable ways: first, the convergence rate of PA is defined in \cite{Fagnani08} as the decay rate of the {\em expected value} of a Lyapunov-like function $d(k)$. Although this stochastic measure captures the average behavior of PA, it offers little guarantee on the decay rate of each realization $(d(k))_{k=0}^\infty$. In contrast, the bounds $\gamma$ on convergence rate of ICHA above are deterministic, providing guarantees on the decay rate of $(V(\mathbf{x}(k)))_{k=0}^\infty$. Second, even if the first difference is disregarded, the bounds of ICHA are still roughly $20$\% better than the convergence rate of PA for a few common graphs. To justify this claim, let $1-\frac{1}{\gamma_{\text{PA}}}$ denote the convergence rate of PA. Since PA requires two real-number transmissions per iteration while ICHA requires only one, to enable a fair comparison we introduce a two-iteration bound $\gamma_{\text{ICHA}}$ for ICHA, defined as $\gamma_{\text{ICHA}}=\frac{\gamma^2}{2\gamma-1}$ so that $1-\frac{1}{\gamma_{\text{ICHA}}}=(1-\frac{1}{\gamma})^2$. Figure~\ref{fig:gammagamma_N} plots the ratio $\frac{\gamma_{\text{ICHA}}}{\gamma_{\text{PA}}}$ versus $N$ for path, cycle, and complete graphs, where $\gamma_{\text{PA}}$ is computed according to \cite{Fagnani08}, while $\gamma_{\text{ICHA}}$ is computed using $\gamma$ in~\ref{enu:specpath}, \ref{enu:speccycl}, and~\ref{enu:genecomp}. Observe that for $N>50$, $\gamma_{\text{ICHA}}$ is $18$\% smaller than $\gamma_{\text{PA}}$ for path and cycle graphs, and $25$\% so for complete graphs. The latter can also be shown analytically: since $\gamma_{\text{PA}}=N-1$ and $\gamma_{\text{ICHA}}=\frac{(\frac{3}{2}N-1)^2}{2(\frac{3}{2}N-1)-1}$, $\lim_{N\rightarrow\infty}\frac{\gamma_{\text{ICHA}}}{\gamma_{\text{PA}}}=\frac{3}{4}$. This justifies the claim. Finally, unlike $\gamma$ and $\gamma_{\text{ICHA}}$, $\gamma_{\text{PA}}$ in general cannot be expressed in a form that explicitly reveals its dependence on the graph invariants. Indeed, it generally can only be computed by numerically finding the spectral radius of an invariant subspace of an $N^2$-by-$N^2$ matrix, which may be prohibitive for large $N$.

\subsection{Practical Version}\label{ssec:pracvers}

The strong convergence properties of ICHA suggest that its greedy behavior may be worthy of emulating. In this subsection, we derive a practical algorithm that closely mimics such behavior.

Reconsider the system \eqref{eq:x=sumcxsumcuijx}, \eqref{eq:xh=sumcxsumc}, \eqref{eq:c=1|N|1|N|}, \eqref{eq:x=y|N|y|N|c} and suppose this system evolves in a discrete event fashion, according to the following description: associated with the system is {\em time}, which is real-valued, nonnegative, and denoted as $t\in[0,\infty)$, where $t=0$ represents the time instant at which the nodes have observed the $y_i$'s but have yet to execute an iteration. In addition, associated with each node $i\in\mathcal{V}$ is an {\em event}, which is scheduled to occur at time $\tau_i\in(0,\infty]$ and is marked by node $i$ initiating an iteration, where $\tau_i=\infty$ means the event will not occur. Each event time $\tau_i$ is a {\em variable}, which is initialized at time $t=0$ to $\tau_i(0)$, is updated only at each iteration $k\in\mathbb{P}$ from $\tau_i(k-1)$ to $\tau_i(k)$, and is no less than $t$ at any time $t$, so that no event is scheduled to occur in the past. Starting from $t=0$, time advances to $t=\min_{i\in\mathcal{V}}\tau_i(0)$, at which an event, marked by node $u(1)\in\operatornamewithlimits{arg\,min}_{i\in\mathcal{V}}\tau_i(0)$ initiating iteration $1$, occurs, during which $\tau_i(1)$ $\forall i\in\mathcal{V}$ are determined. Time then advances to $t=\min_{i\in\mathcal{V}}\tau_i(1)$, at which a subsequent event, marked by node $u(2)\in\operatornamewithlimits{arg\,min}_{i\in\mathcal{V}}\tau_i(1)$ initiating iteration $2$, occurs, during which $\tau_i(2)$ $\forall i\in\mathcal{V}$ are determined. In the same way, time continues to advance toward infinity, while events continue to occur one after another, except if $\tau_i(k)=\infty$ $\forall i\in\mathcal{V}$ for some $k\in\mathbb{N}$, for which the system terminates.

Having described how the system evolves, we now specify how $\tau_i(k)$ $\forall k\in\mathbb{N}$ $\forall i\in\mathcal{V}$ are recursively determined. First, consider the time instant $t=0$, at which $\tau_i(0)$ $\forall i\in\mathcal{V}$ need to be determined. To behave greedily, nodes with the maximum $\Delta V_i(\mathbf{x}(0))$'s should have the minimum $\tau_i(0)$'s. This may be accomplished by letting
\begin{align}
\tau_i(0)=\Phi(\Delta V_i(\mathbf{x}(0))),\quad\forall i\in\mathcal{V},\label{eq:t=PDV}
\end{align}
where $\Phi:[0,\infty)\rightarrow(0,\infty]$ is a continuous and strictly decreasing function satisfying $\lim_{v\rightarrow0}\Phi(v)=\infty$ and $\Phi(0)=\infty$. Although, mathematically, \eqref{eq:t=PDV} ensures that $V(\mathbf{x}(0))$ drops maximally to $V(\mathbf{x}(1))$, in reality it is possible that multiple nodes have the same minimum $\tau_i(0)$'s, leading to wireless collisions. To address this issue, we insert a little randomness into \eqref{eq:t=PDV}, rewriting it as
\begin{align}
\tau_i(0)=\Phi(\Delta V_i(\mathbf{x}(0)))+\varepsilon(\Delta V_i(\mathbf{x}(0)))\cdot\operatorname{rand}(),\quad\forall i\in\mathcal{V},\label{eq:t0=PDVeDVr}
\end{align}
where $\varepsilon:[0,\infty)\rightarrow(0,\infty)$ is a continuous function meant to take on small positive values and each call to $\operatorname{rand}()$ returns a uniformly distributed random number in $(0,1)$. With \eqref{eq:t0=PDVeDVr}, with high probability iteration $1$ is initiated by a node $i$ with the maximum, or a near-maximum, $\Delta V_i(\mathbf{x}(0))$.

Next, pick any $k\in\mathbb{P}$ and consider the time instant $t=\min_{i\in\mathcal{V}}\tau_i(k-1)$, at which node $u(k)\in\operatornamewithlimits{arg\,min}_{i\in\mathcal{V}}\tau_i(k-1)$ initiates iteration $k$, during which $\tau_i(k)$ $\forall i\in\mathcal{V}$ need to be determined. Again, to be greedy, nodes with the maximum $\Delta V_i(\mathbf{x}(k))$'s should have the minimum $\tau_i(k)$'s. At first glance, this may be approximately accomplished following ideas from \eqref{eq:t0=PDVeDVr}, i.e., by letting
\begin{align}
\tau_i(k)=\Phi(\Delta V_i(\mathbf{x}(k)))+\varepsilon(\Delta V_i(\mathbf{x}(k)))\cdot\operatorname{rand}(),\quad\forall i\in\mathcal{V}.\label{eq:t=PDVeDVr}
\end{align}
However, with \eqref{eq:t=PDVeDVr}, it is possible that $\tau_i(k)$ turns out to be smaller than $t$, causing an event to be scheduled in the past. Moreover, nodes who are two or more hops away from node $u(k)$ are unaware of the ongoing iteration $k$ and, thus, are unable to perform an update. Fortunately, these issues may be overcome by slightly modifying \eqref{eq:t=PDVeDVr} as follows:
\begin{align}
\tau_i(k)=\begin{cases}\max\{\Phi(\Delta V_i(\mathbf{x}(k))),t\}+\varepsilon(\Delta V_i(\mathbf{x}(k)))\cdot\operatorname{rand}(), & \text{if $i\in\mathcal{N}_{u(k)}\cup\{u(k)\}$},\\ \tau_i(k-1), & \text{otherwise},\end{cases}\quad\forall i\in\mathcal{V}.\label{eq:t=maxPDVteDVrt}
\end{align}
Using \eqref{eq:t0=PDVeDVr} and \eqref{eq:t=maxPDVteDVrt} and by induction on $k'\in\mathbb{P}$, it can be shown that $\tau_i(k')$ satisfies
\begin{align*}
\max\{\Phi(\Delta V_i(\mathbf{x}(k'))),t'\}\le\tau_i(k')\le\max\{\Phi(\Delta V_i(\mathbf{x}(k'))),t'\}+\varepsilon(\Delta V_i(\mathbf{x}(k'))),\quad\forall k'\in\mathbb{P},\;\forall i\in\mathcal{V},
\end{align*}
where $t'=\min_{j\in\mathcal{V}}\tau_j(k'-1)$. Hence, with \eqref{eq:t=maxPDVteDVrt}, it is highly probable that iteration $k+1$ is initiated by a node $i$ with the maximum or a near-maximum $\Delta V_i(\mathbf{x}(k))$. It follows that with \eqref{eq:t0=PDVeDVr} and \eqref{eq:t=maxPDVteDVrt}, the nodes closely mimic the greedy behavior of ICHA. Note that \eqref{eq:t0=PDVeDVr} and \eqref{eq:t=maxPDVteDVrt} represent a {\em feedback iteration controller}, which uses architecture~\ref{enu:useDVfeedcontwheninit} and follows the spirit of~\ref{enu:largDVsooninit} (since $\Phi$ is strictly decreasing and $\varepsilon$ is small) and~\ref{enu:whenDV0refrinit} (since $\Phi(0)=\infty$). Also, $\Phi$ and $\varepsilon$ represent the {\em controller parameters}, which may be selected based on practical wireless networking considerations (e.g., all else being equal, $\Phi(v)=\frac{1}{v}$ and $\varepsilon(v)=0.001$ yield faster convergence time than $\Phi(v)=\frac{10}{v}$ and $\varepsilon(v)=0.01$ but higher collision probability).

The above description defines a discrete event system, which can be realized via a distributed asynchronous algorithm, referred to as {\em Controlled Hopwise Averaging} (CHA) and stated as follows:

\begin{algorithm}[Controlled Hopwise Averaging]\label{alg:CHA}
\begin{algorithminit}{}
\item Let time $t=0$.
\item Each node $i\in\mathcal{V}$ transmits $|\mathcal{N}_i|$ and $y_i$ to every node $j\in\mathcal{N}_i$.
\item Each node $i\in\mathcal{V}$ creates variables $x_{ij}\in\mathbb{R}$ $\forall j\in\mathcal{N}_i$, $\hat{x}_i\in\mathbb{R}$, $\Delta V_i\in[0,\infty)$, and $\tau_i\in(0,\infty]$ and initializes them sequentially:
\\$x_{ij}\leftarrow\frac{\frac{y_i}{|\mathcal{N}_i|}+\frac{y_j}{|\mathcal{N}_j|}}{c_{\{i,j\}}},\quad\forall j\in\mathcal{N}_i,$
\\$\hat{x}_i\leftarrow\frac{\sum_{j\in\mathcal{N}_i}c_{\{i,j\}}x_{ij}}{\sum_{j\in\mathcal{N}_i}c_{\{i,j\}}},$
\\$\Delta V_i\leftarrow\sum_{j\in\mathcal{N}_i}c_{\{i,j\}}(x_{ij}-\hat{x}_i)^2,$
\\$\tau_i\leftarrow\Phi(\Delta V_i)+\varepsilon(\Delta V_i)\cdot\operatorname{rand}().$
\end{algorithminit}
\begin{algorithmoper}{At each iteration:}
\item Let $t=\min_{j\in\mathcal{V}}\tau_j$ and $i\in\operatornamewithlimits{arg\,min}_{j\in\mathcal{V}}\tau_j$.
\item Node $i$ updates $x_{ij}$ $\forall j\in\mathcal{N}_i$, $\Delta V_i$, and $\tau_i$ sequentially:
\\$x_{ij}\leftarrow\hat{x}_i,\quad\forall j\in\mathcal{N}_i,$
\\$\Delta V_i\leftarrow0,$
\\$\tau_i\leftarrow\infty.$
\item Node $i$ transmits $\hat{x}_i$ to every node $j\in\mathcal{N}_i$.
\item Each node $j\in\mathcal{N}_i$ updates $x_{ji}$, $\hat{x}_j$, $\Delta V_j$, and $\tau_j$ sequentially:
\\$x_{ji}\leftarrow\hat{x}_i,$
\\$\hat{x}_j\leftarrow\frac{\sum_{\ell\in\mathcal{N}_j}c_{\{j,\ell\}}x_{j\ell}}{\sum_{\ell\in\mathcal{N}_j}c_{\{j,\ell\}}},$
\\$\Delta V_j\leftarrow\sum_{\ell\in\mathcal{N}_j}c_{\{j,\ell\}}(x_{j\ell}-\hat{x}_j)^2,$
\\$\tau_j\leftarrow\max\{\Phi(\Delta V_j),t\}+\varepsilon(\Delta V_j)\cdot\operatorname{rand}().$
\end{algorithmoper}
\end{algorithm}

Algorithm~\ref{alg:CHA}, or CHA, is similar to ICHA in Algorithm~\ref{alg:ICHA} except that each node $i$ maintains an additional variable $\tau_i$, in Steps~3, 5, and~7, and that each iteration is initiated, in a discrete event fashion, by a node $i$ having the minimum $\tau_i$, in Step~4. Note that ``$\tau_i\leftarrow\infty$'' in Step~5 is due to ``$\Delta V_i\leftarrow0$'' and to $\Phi(0)=\infty$. Moreover, every step of CHA is implementable in a fully decentralized manner, making it a practical algorithm.

To analyze the behavior of CHA, recall that $\varepsilon$ is meant to take on small positive values, creating just a little randomness so that the probability of wireless collisions is zero. For the purpose of analysis, we turn this feature off (i.e., set $\varepsilon(v)=0$ $\forall v\in[0,\infty)$) and let the symbol ``$\in$'' in Step~4 take care of the randomness (i.e., randomly pick an element $i$ from the set $\operatornamewithlimits{arg\,min}_{j\in\mathcal{V}}\tau_j$ whenever it has multiple elements). We also allow $\Phi$ to be arbitrary (but satisfy the conditions stated when it was introduced). With this setup, the following convergence properties of CHA can be established:

\begin{theorem}\label{thm:CHAexpconv}
Theorems~\ref{thm:ICHAexpconvgene} and~\ref{thm:ICHAexpconvspec}, intended for ICHA described in Algorithm~\ref{alg:ICHA}, hold verbatim for CHA described in Algorithm~\ref{alg:CHA} with any $\Phi$ and with $\varepsilon$ satisfying $\varepsilon(v)=0$ $\forall v\in[0,\infty)$. In addition, $\lim_{k\rightarrow\infty}t(k)=\infty$ and $V(\mathbf{x}(k))\le(\gamma-1)\Phi^{-1}(t(k))$ $\forall k\in\mathbb{P}$, where $t(0)=0$ and $t(k)$ is the time instant at which iteration $k$ occurs.
\end{theorem}

\begin{proof}
See Appendix~\ref{ssec:proofthmCHAexpconv}.
\end{proof}

Theorem~\ref{thm:CHAexpconv} characterizes the convergence of CHA in two senses: {\em iteration} and {\em time}. Iteration-wise, it says that CHA converges exponentially and shares the same bounds $\gamma$ on convergence rate as ICHA, regardless of $\Phi$. This result suggests that CHA does closely emulate ICHA. Time-wise, the theorem says that CHA converges asymptotically and perhaps exponentially, depending on $\Phi$. For example, $\Phi(v)=\frac{1}{v}$ does not guarantee exponential convergence in time (since $\Phi^{-1}(v)=\frac{1}{v}$), but $\Phi(v)=W(\frac{1}{v})$, where $W$ is the Lambert W function, does (since $\Phi^{-1}(v)=\frac{1}{v}e^{-v}$). Therefore, the controller parameter $\Phi$ may be used to shape the temporal convergence of CHA.

\begin{remark}\label{rem:clocoffs}
CHA has a limitation: it assumes no clock offsets among the nodes. Note, however, that although such offsets would cause CHA to deviate from its designed behavior, they would not render it ``inoperable,'' i.e., $V(\mathbf{x}(k))$ would still strictly decrease after every iteration $k$, and the conservation \eqref{eq:sumcx=sumcx} would still hold, so that the $x_{\{i,j\}}(k)$'s and $\hat{x}_i(k)$'s would still approach $x^*$.
\end{remark}

\section{Performance Comparison}\label{sec:perfcomp}

In this section, we compare the performance of RHA and CHA with that of Pairwise Averaging (PA) \cite{Tsitsiklis84}, Consensus Propagation (CP) \cite{Moallemi06}, Algorithm A2 (A2) of \cite{Mehyar07}, and Distributed Random Grouping (DRG) \cite{ChenJY06} via extensive simulation on multi-hop wireless networks modeled by random geometric graphs. For completeness, PA, CP, A2, and DRG are stated below, in which $\mathcal{E}'=\{(i,j)\in\mathcal{V}\times\mathcal{V}:\{i,j\}\in\mathcal{E}\}$ denotes the set of $2L$ directed links:

\begin{algorithm}[Pairwise Averaging \cite{Tsitsiklis84}]\label{alg:PA}
\begin{algorithminit}{}
\item Each node $i\in\mathcal{V}$ creates a variable $\hat{x}_i\in\mathbb{R}$ and initializes it: $\hat{x}_i\leftarrow y_i$.
\end{algorithminit}
\begin{algorithmoper}{At each iteration:}
\item A link, say, link $\{i,j\}$, is selected randomly and equiprobably out of the set $\mathcal{E}$ of $L$ links. Node $i$ transmits $\hat{x}_i$ to node $j$. Node $j$ updates $\hat{x}_j$: $\hat{x}_j\leftarrow\frac{\hat{x}_i+\hat{x}_j}{2}$. Node $j$ transmits $\hat{x}_j$ to node $i$. Node $i$ updates $\hat{x}_i$: $\hat{x}_i\leftarrow\hat{x}_j$.
\end{algorithmoper}
\end{algorithm}

\begin{algorithm}[Consensus Propagation \cite{Moallemi06}]\label{alg:CP}
\begin{algorithminit}{}
\item Each node $i\in\mathcal{V}$ creates variables $K_{ji}\ge0$ $\forall j\in\mathcal{N}_i$, $\mu_{ji}\in\mathbb{R}$ $\forall j\in\mathcal{N}_i$, and $\hat{x}_i\in\mathbb{R}$ and initializes them sequentially: $K_{ji}\leftarrow0$ $\forall j\in\mathcal{N}_i$, $\mu_{ji}\leftarrow0$ $\forall j\in\mathcal{N}_i$, $\hat{x}_i\leftarrow y_i$.
\end{algorithminit}
\begin{algorithmoper}{At each iteration:}
\item A directed link, say, link $(i,j)$, is selected randomly and equiprobably out of the set $\mathcal{E}'$ of $2L$ directed links. Node $i$ transmits $F_{ij}\triangleq\frac{1+\sum_{\ell\in\mathcal{N}_i,\ell\neq j}K_{\ell i}}{1+\frac{1}{\beta}(1+\sum_{\ell\in\mathcal{N}_i,\ell\neq j}K_{\ell i})}$ and $G_{ij}\triangleq\frac{y_i+\sum_{\ell\in\mathcal{N}_i,\ell\neq j}K_{\ell i}\mu_{\ell i}}{1+\sum_{\ell\in\mathcal{N}_i,\ell\neq j}K_{\ell i}}$ to node $j$. Node $j$ updates $K_{ij}$, $\mu_{ij}$, and $\hat{x}_j$ sequentially: $K_{ij}\leftarrow F_{ij}$, $\mu_{ij}\leftarrow G_{ij}$, $\hat{x}_j\leftarrow\frac{y_j+\sum_{\ell\in\mathcal{N}_j}K_{\ell j}\mu_{\ell j}}{1+\sum_{\ell\in\mathcal{N}_j}K_{\ell j}}$.
\end{algorithmoper}
\end{algorithm}

\begin{algorithm}[Algorithm A2 \cite{Mehyar07}]\label{alg:A2}
\begin{algorithminit}{}
\item Each node $i\in\mathcal{V}$ creates variables $\delta_{ij}\in\mathbb{R}$ $\forall j\in\mathcal{N}_i$ and $\hat{x}_i\in\mathbb{R}$ and initializes them sequentially: $\delta_{ij}\leftarrow0$ $\forall j\in\mathcal{N}_i$, $\hat{x}_i\leftarrow y_i$.
\end{algorithminit}
\begin{algorithmoper}{At each iteration:}
\item A directed link, say, link $(i,j)$, is selected randomly and equiprobably out of the set $\mathcal{E}'$ of $2L$ directed links. Node $i$ transmits $\hat{x}_i$ to node $j$. Node $j$ updates $\delta_{ji}$: $\delta_{ji}\leftarrow\delta_{ji}+\phi(\hat{x}_i-\hat{x}_j)$. Node $j$ transmits $\phi(\hat{x}_i-\hat{x}_j)$ to node $i$. Node $i$ updates $\delta_{ij}$: $\delta_{ij}\leftarrow\delta_{ij}-\phi(\hat{x}_i-\hat{x}_j)$. Each node $\ell\in\mathcal{V}$ updates $\hat{x}_\ell$: $\hat{x}_\ell\leftarrow\hat{x}_\ell+\frac{\gamma}{|\mathcal{N}_\ell|+1}((\sum_{m\in\mathcal{N}_\ell}\delta_{\ell m})+y_\ell-\hat{x}_\ell)$.
\end{algorithmoper}
\end{algorithm}

\begin{algorithm}[Distributed Random Grouping \cite{ChenJY06}]\label{alg:DRG}
\begin{algorithminit}{}
\item Each node $i\in\mathcal{V}$ creates a variable $\hat{x}_i\in\mathbb{R}$ and initializes it: $\hat{x}_i\leftarrow y_i$.
\end{algorithminit}
\begin{algorithmoper}{At each iteration:}
\item A node, say, node $i$, is selected randomly and equiprobably out of the set $\mathcal{V}$ of $N$ nodes. Node $i$ transmits a message to every node $j\in\mathcal{N}_i$, requesting their $\hat{x}_j$'s. Each node $j\in\mathcal{N}_i$ transmits $\hat{x}_j$ to node $i$. Node $i$ updates $\hat{x}_i$: $\hat{x}_i\leftarrow\frac{\sum_{j\in\{i\}\cup\mathcal{N}_i}\hat{x}_j}{|\mathcal{N}_i|+1}$. Node $i$ transmits $\hat{x}_i$ to every node $j\in\mathcal{N}_i$. Each node $j\in\mathcal{N}_i$ updates $\hat{x}_j$: $\hat{x}_j\leftarrow\hat{x}_i$.
\end{algorithmoper}
\end{algorithm}

Note that RHA and CHA require $2N$ real-number transmissions as initialization overhead, whereas PA, CP, A2, and DRG require none. However, PA, CP, and A2 require two real-number transmissions per iteration and DRG requires $|\mathcal{N}_i|+1$ (where $i$ is the node that leads an iteration), whereas RHA and CHA require only one. Also note that CP has a parameter $\beta\in(0,\infty]$ and A2 has two parameters $\gamma\in(0,1)$ and $\phi\in(0,\frac{1}{2})$. Moreover, PA and DRG are assumed to be free of overlapping iterations, i.e., deficiency~\ref{enu:overiter}.

\begin{figure}[tb]
\centering\includegraphics{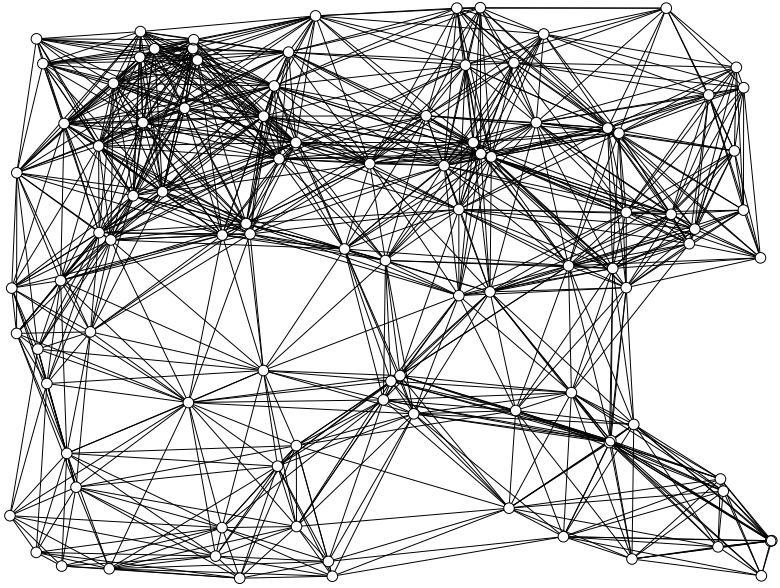}
\caption{A $100$-node, $1000$-link multi-hop wireless network.}
\label{fig:net}
\end{figure}

To compare the performance of these algorithms, two sets of simulation are carried out. The first set corresponds to a single scenario of a multi-hop wireless network with $N=100$ nodes, where each node $i$ observes $y_i\in(0,1)$ and has, on average, $\frac{2L}{N}=20$ one-hop neighbors, as shown in Figure~\ref{fig:net}. The second set corresponds to multi-hop wireless networks modeled by random geometric graphs, with the number of nodes varying from $N=100$ to $N=500$, and the average number of neighbors varying from $\frac{2L}{N}=10$ to $\frac{2L}{N}=60$. For each $N$ and $\frac{2L}{N}$, we generate $50$ scenarios. For each scenario, we randomly and uniformly place $N$ nodes in the unit square $(0,1)\times(0,1)$, gradually increase the one-hop radius until there are $L$ links (or $\frac{2L}{N}$ neighbors on average), randomly and uniformly generate the $y_i$'s in $(0,1)$, and repeat this process if the resulting network is not connected. We then simulate PA, CP, A2, DRG, RHA, and CHA until $3N^2$ real-number transmissions have occurred (i.e., three times of what flooding needs), record the number of real-number transmissions needed to converge (including initialization overhead, if any), and assume that this number is $3N^2$ if an algorithm fails to converge after $3N^2$. For both sets of simulation, we let the convergence criterion be $|\hat{x}_i-x^*|\le0.005$ $\forall i\in\mathcal{V}$ and the parameters be $\beta=10^6$ for CP (obtained after some tuning), $\gamma=0.3$ and $\phi=0.49$ for A2 (ditto), and $\Phi(v)=\frac{1}{v}$ and $\varepsilon(v)=0.001$ for CHA.

\begin{figure}[tb]
\centering\includegraphics{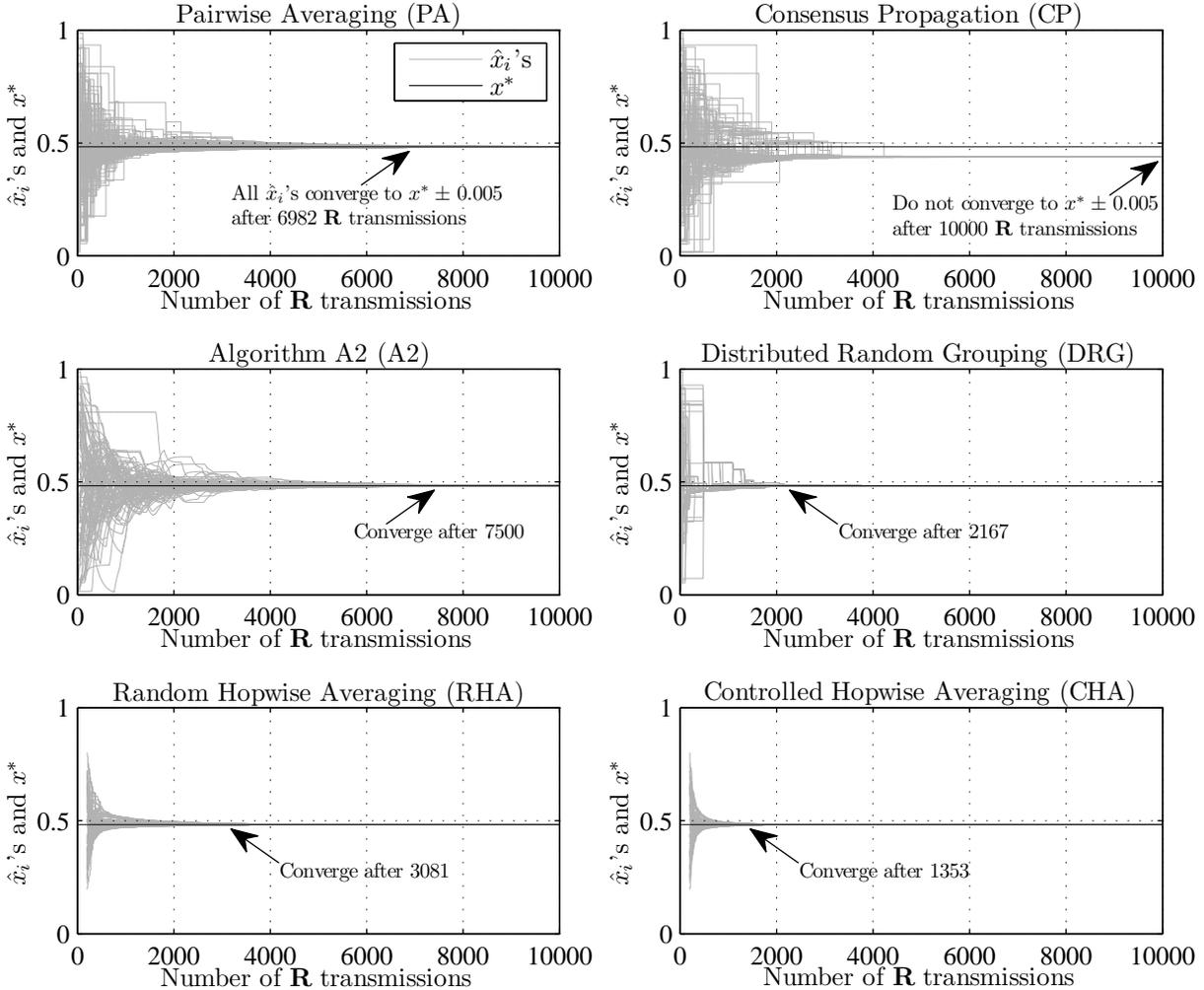}
\caption{Convergence of the estimates $\hat{x}_i(k)$'s to the unknown average $x^*$ under PA, CP, A2, DRG, RHA, and CHA for the network in Figure~\ref{fig:net}.}
\label{fig:xhat_nrt}
\end{figure}

Results from the first set of simulation are shown in Figure~\ref{fig:xhat_nrt}. Observe that PA and A2 have roughly the same performance, requiring approximately $7,000$ real-number transmissions to converge. In contrast, CP fails to converge after $10,000$ transmissions, although it does achieve a consensus. On the other hand, DRG is found to be quite efficient, needing only approximately $2,100$ transmissions for convergence. Note that RHA outperforms PA, CP, and A2, but not DRG, while CHA is the most efficient, requiring only roughly $1,300$ transmissions to converge.

\begin{figure}[tb]
\centering\includegraphics{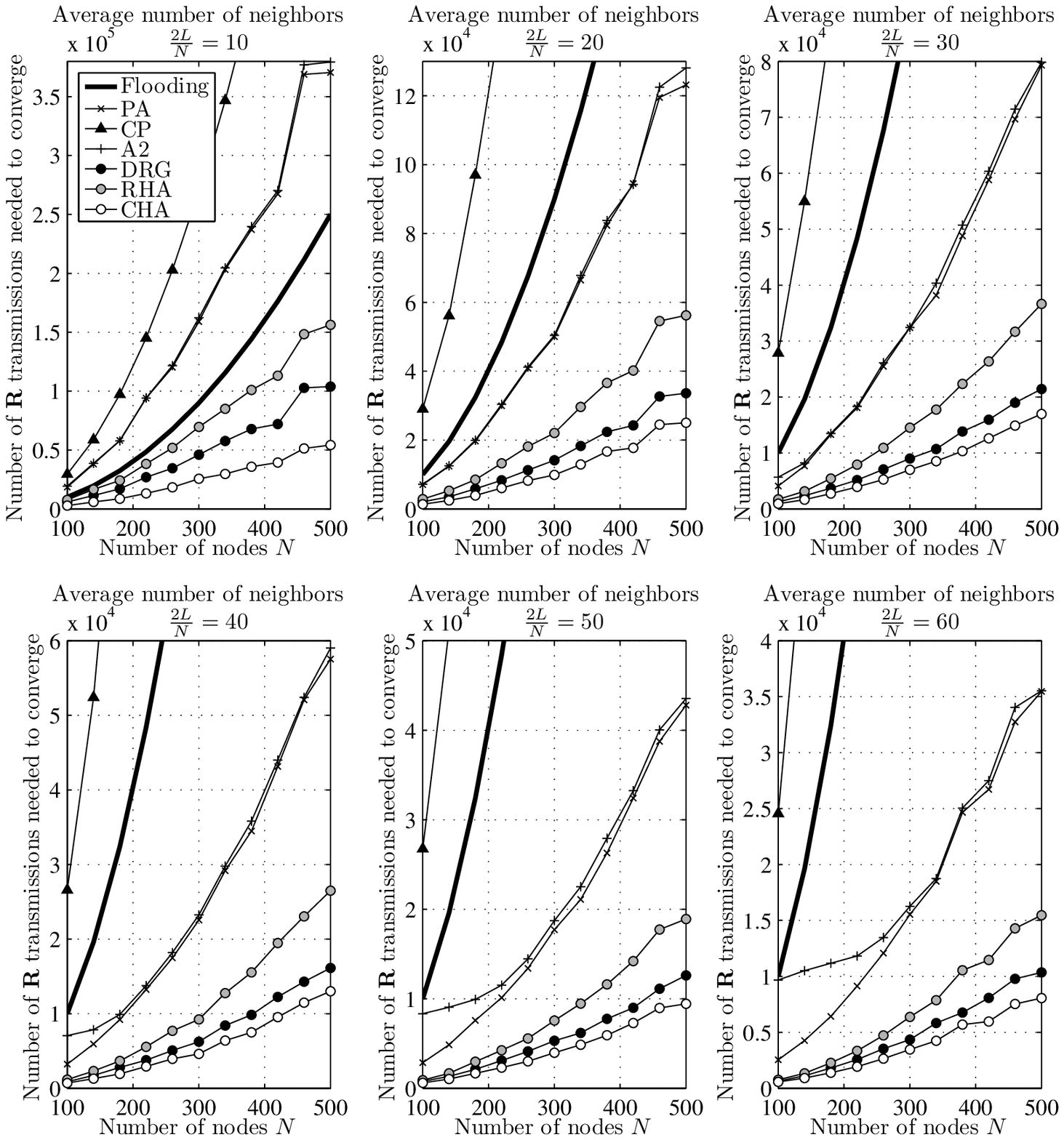}
\caption{Bandwidth/energy efficiency of flooding, PA, CP, A2, DRG, RHA, and CHA on random geometric networks with varying number of nodes $N$ and average number of neighbors $\frac{2L}{N}$.}
\label{fig:nrtc_N}
\end{figure}

Results from the second set of simulation are shown in Figure~\ref{fig:nrtc_N}, where the number of real-number transmissions needed to converge, averaged over $50$ scenarios, is plotted as a function of the number of nodes $N$ and the average number of neighbors $\frac{2L}{N}$. Also included in the figure, as a baseline for comparison, is the performance of flooding (i.e., $N^2$). Observe that regardless of $N$ and $\frac{2L}{N}$, CP has the worst bandwidth/energy efficiency, followed by PA and A2. In addition, DRG, RHA, and CHA are all fairly efficient, with CHA again having the best efficiency. In particular, CHA is at least $20$\% more efficient than DRG, and around $50$\% more so when the network is sparsely connected, at $\frac{2L}{N}=10$. Notice that the performance of DRG is achieved under the assumption that overlapping iterations cannot occur, a condition that CHA does not require. Finally, the significant difference in efficiency between RHA and CHA demonstrates the benefit of incorporating greedy, decentralized, feedback iteration control.

\section{Conclusion}\label{sec:concl}

In this paper, we have shown that the existing distributed averaging schemes have a few drawbacks, which hurt their bandwidth/energy efficiency. Motivated by this, we have devised RHA, an asynchronous algorithm that exploits the broadcast nature of wireless medium, achieves almost sure asymptotic convergence, and overcomes all but one of the drawbacks. To deal with the remaining drawback, on lack of control, we have introduced a new way to apply Lyapunov stability theory, namely, the concept of greedy, decentralized, feedback iteration control. Based on this concept, we have developed ICHA and CHA, established bounds on their exponential convergence rates, and shown that CHA is practical and capable of closely mimicking the behavior of ICHA. Finally, we have shown via extensive simulation that CHA is substantially more bandwidth/energy efficient than several existing schemes.

Several extensions of this work are possible, including design and analysis of ``controlled'' distributed averaging algorithms that are applicable to more general wireless networks (e.g., with directed links, time-varying topologies, and dynamic observations) and more realistic communication channels (e.g., with random delays, packet losses, and quantization effects), and that take into account MAC/PHY layer design issues (e.g., retransmission and backoff strategies).

\appendix
\section{Appendix}\label{sec:app}

\subsection{Proof of Theorem~\ref{thm:ICHAexpconvgene}}\label{ssec:proofthmICHAexpconvgene}

To prove Theorem~\ref{thm:ICHAexpconvgene}, we first prove the following lemma:

\begin{lemma}\label{lem:V<=gmaxDVgene}
$V(\mathbf{x}(k))\le\gamma\max_{i\in\mathcal{V}}\Delta V_i(\mathbf{x}(k))$ $\forall k\in\mathbb{N}$, where $\gamma$ is as in \eqref{eq:g=N2aN2b3N1DD12N}.
\end{lemma}

\begin{proof}
Let $k\in\mathbb{N}$. Notice from \eqref{eq:c=1|N|1|N|} that $\sum_{i\in\mathcal{V}}b_i=N$ and from \eqref{eq:x*=1Nsumy}, \eqref{eq:xh=sumcxsumc}, \eqref{eq:sumcx=sumcx}, and \eqref{eq:sumcxsumc=sumyN} that $\sum_{i\in\mathcal{V}}b_i\hat{x}_i(k)=\sum_{i\in\mathcal{V}}b_ix^*$. Thus, $\sum_{i\in\mathcal{V}}\sum_{j\in\mathcal{V}}b_ib_j(\hat{x}_i(k)-\hat{x}_j(k))^2=\sum_{j\in\mathcal{V}}b_j\sum_{i\in\mathcal{V}}b_i(\hat{x}_i(k)-x^*)^2+\sum_{i\in\mathcal{V}}b_i\sum_{j\in\mathcal{V}}b_j(\hat{x}_j(k)-x^*)^2-2\sum_{i\in\mathcal{V}}b_i(\hat{x}_i(k)-x^*)\sum_{j\in\mathcal{V}}b_j(\hat{x}_j(k)-x^*)=2N\sum_{i\in\mathcal{V}}b_i(\hat{x}_i(k)-x^*)^2$. It follows from \eqref{eq:V=sumcxx*2}, \eqref{eq:DV=sumcxxh2}, and \eqref{eq:xh=sumcxsumc} that
\begin{align}
&V(\mathbf{x}(k))=\frac{1}{2}\sum_{i\in\mathcal{V}}\sum_{j\in\mathcal{N}_i}c_{\{i,j\}}(x_{\{i,j\}}(k)-\hat{x}_i(k))^2+\frac{1}{2}\sum_{i\in\mathcal{V}}\sum_{j\in\mathcal{N}_i}c_{\{i,j\}}(\hat{x}_i(k)-x^*)^2\nonumber\displaybreak[0]\\
&\quad+\sum_{i\in\mathcal{V}}(\hat{x}_i(k)-x^*)\sum_{j\in\mathcal{N}_i}c_{\{i,j\}}(x_{\{i,j\}}(k)-\hat{x}_i(k))=\frac{1}{2}\sum_{i\in\mathcal{V}}\Delta V_i(\mathbf{x}(k))+\sum_{i\in\mathcal{V}}b_i(\hat{x}_i(k)-x^*)^2\label{eq:V=12sumDVsumbxhx*2}\displaybreak[0]\\
&\quad\le\tfrac{N}{2}\max_{i\in\mathcal{V}}\Delta V_i(\mathbf{x}(k))+\tfrac{\sum_{i\in\mathcal{V}}\sum_{j\in\mathcal{V}}b_ib_j(\hat{x}_i(k)-\hat{x}_j(k))^2}{2N}\nonumber\displaybreak[0]\\
&\quad=\tfrac{N}{2}\max_{i\in\mathcal{V}}\Delta V_i(\mathbf{x}(k))+\tfrac{\sum_{i\in\mathcal{V}}\sum_{j\in\mathcal{N}_i}b_ib_j(\hat{x}_i(k)-\hat{x}_j(k))^2}{2N}+\tfrac{\sum_{i\in\mathcal{V}}\sum_{j\in\mathcal{V}-\mathcal{N}_i-\{i\}}b_ib_j(\hat{x}_i(k)-\hat{x}_j(k))^2}{2N}.\label{eq:N2maxDVsumsumbbxhxh22Nsumsumbbxhxh22N}
\end{align}
Note from \eqref{eq:DV=sumcxxh2} that $N\max_{i\in\mathcal{V}}\Delta V_i(\mathbf{x}(k))\ge\sum_{i\in\mathcal{V}}b_i\Delta V_i(\mathbf{x}(k))=\sum_{\{i,j\}\in\mathcal{E}}b_ic_{\{i,j\}}(\hat{x}_i(k)-x_{\{i,j\}}(k))^2+b_jc_{\{i,j\}}(\hat{x}_j(k)-x_{\{i,j\}}(k))^2\ge\sum_{\{i,j\}\in\mathcal{E}}\frac{b_ib_jc_{\{i,j\}}}{b_i+b_j}(\hat{x}_i(k)-\hat{x}_j(k))^2$. Hence,
\begin{align}
\sum_{i\in\mathcal{V}}\sum_{j\in\mathcal{N}_i}b_ib_j(\hat{x}_i(k)-\hat{x}_j(k))^2\le2\alpha N\max_{i\in\mathcal{V}}\Delta V_i(\mathbf{x}(k)).\label{eq:sumsumbbxhxh2<=2aNmaxDV}
\end{align}
Next, it can be shown via \eqref{eq:DV=sumcxxh2} that $\forall i\in\mathcal{V}$ with $|\mathcal{N}_i|\ge2$, $\forall j,\ell\in\mathcal{N}_i$ with $j\neq\ell$, $c_{\{i,j\}}c_{\{i,\ell\}}(x_{\{i,j\}}(k)-x_{\{i,\ell\}}(k))^2\le(c_{\{i,j\}}+c_{\{i,\ell\}})(c_{\{i,j\}}(x_{\{i,j\}}(k)-\hat{x}_i(k))^2+c_{\{i,\ell\}}(x_{\{i,\ell\}}(k)-\hat{x}_i(k))^2)\le(c_{\{i,j\}}+c_{\{i,\ell\}})\Delta V_i(\mathbf{x}(k))$, implying that $|x_{\{i,j\}}(k)-x_{\{i,\ell\}}(k)|\le\bigl(\max_{p\in\mathcal{V}}\Delta V_p(\mathbf{x}(k))(\frac{1}{c_{\{i,j\}}}+\frac{1}{c_{\{i,\ell\}}})\bigr)^{\frac{1}{2}}$. In addition, $\forall i\in\mathcal{V}$, $\forall j\in\mathcal{N}_i$, $|\hat{x}_i(k)-x_{\{i,j\}}(k)|\le\bigl(\frac{\max_{p\in\mathcal{V}}\Delta V_p(\mathbf{x}(k))}{c_{\{i,j\}}}\bigr)^{\frac{1}{2}}$ because of \eqref{eq:DV=sumcxxh2}. For any $i,j\in\mathcal{V}$ with $i\neq j$, let the sequence $(a_1,a_2,\ldots,a_{m_{ij}})$ represent a shortest path from node $i$ to node $j$, where $a_1=i$, $a_{m_{ij}}=j$, $\{a_\ell,a_{\ell+1}\}\in\mathcal{E}$ $\forall \ell\in\{1,2,\ldots,m_{ij}-1\}$, and $2\le m_{ij}\le D+1$. Then, it follows from \eqref{eq:c=1|N|1|N|}, the triangle inequality, and the root-mean square-arithmetic mean-geometric mean inequality that $|\hat{x}_i(k)-\hat{x}_j(k)|\le\bigl(\max_{p\in\mathcal{V}}\Delta V_p(\mathbf{x}(k))\bigr)^{\frac{1}{2}}\Bigl(\bigl(\frac{|\mathcal{N}_{a_1}|\cdot|\mathcal{N}_{a_2}|}{|\mathcal{N}_{a_1}|+|\mathcal{N}_{a_2}|}\bigr)^{\frac{1}{2}}+\sum_{\ell=2}^{m_{ij}-1}\bigl(\frac{|\mathcal{N}_{a_{\ell-1}}|\cdot|\mathcal{N}_{a_\ell}|}{|\mathcal{N}_{a_{\ell-1}}|+|\mathcal{N}_{a_\ell}|}+\frac{|\mathcal{N}_{a_\ell}|\cdot|\mathcal{N}_{a_{\ell+1}}|}{|\mathcal{N}_{a_\ell}|+|\mathcal{N}_{a_{\ell+1}}|}\bigr)^{\frac{1}{2}}+\bigl(\frac{|\mathcal{N}_{a_{m_{ij}-1}}|\cdot|\mathcal{N}_{a_{m_{ij}}}|}{|\mathcal{N}_{a_{m_{ij}-1}}|+|\mathcal{N}_{a_{m_{ij}}}|}\bigr)^{\frac{1}{2}}\Bigr)\le\bigl(m_{ij}\max_{p\in\mathcal{V}}\Delta V_p(\mathbf{x}(k))\bigr)^{\frac{1}{2}}\Bigl(\frac{|\mathcal{N}_{a_1}|+|\mathcal{N}_{a_2}|}{4}+\sum_{\ell=2}^{m_{ij}-1}\bigl(\frac{|\mathcal{N}_{a_{\ell-1}}|+|\mathcal{N}_{a_\ell}|}{4}+\frac{|\mathcal{N}_{a_\ell}|+|\mathcal{N}_{a_{\ell+1}}|}{4}\bigr)+\frac{|\mathcal{N}_{a_{m_{ij}-1}}|+|\mathcal{N}_{a_{m_{ij}}}|}{4}\Bigr)^{\frac{1}{2}}\le\bigl(m_{ij}\max_{p\in\mathcal{V}}\Delta V_p(\mathbf{x}(k))\sum_{\ell=1}^{m_{ij}}|\mathcal{N}_{a_\ell}|\bigr)^{\frac{1}{2}}$. Next, we show that $\forall i,j\in\mathcal{V}$ with $i\neq j$, each node $\ell\in\mathcal{V}-\{a_1,a_2,\ldots,a_{m_{ij}}\}$ has at most $3$ one-hop neighbors in $\{a_1,a_2,\ldots,a_{m_{ij}}\}$. Clearly, this statement is true for $m_{ij}\le3$. For $m_{ij}\ge4$, assume to the contrary that $\exists\ell\in\mathcal{V}-\{a_1,a_2,\ldots,a_{m_{ij}}\}$ such that $\mathcal{N}_\ell\cap\{a_1,a_2,\ldots,a_{m_{ij}}\}=\{a_{i_1},a_{i_2},\ldots,a_{i_n}\}$ for some $1\le i_1<i_2<\cdots<i_n\le m_{ij}$ and $n\ge4$. Then, $(a_1,\ldots,a_{i_1},\ell,a_{i_n},\ldots,a_{m_{ij}})$ is a path shorter than the shortest path $(a_1,a_2,\ldots,a_{m_{ij}})$, which is a contradiction. Therefore, the statement is true. Consequently, $\sum_{\ell=1}^{m_{ij}}|\mathcal{N}_{a_\ell}|\le3(N-m_{ij})+2(m_{ij}-1)=3N-m_{ij}-2$. It follows that $\forall i,j\in\mathcal{V}$ with $i\neq j$, $(\hat{x}_i(k)-\hat{x}_j(k))^2\le m_{ij}(3N-m_{ij}-2)\max_{p\in\mathcal{V}}\Delta V_p(\mathbf{x}(k))$. Since $m_{ij}\le D+1\le N$, $(\hat{x}_i(k)-\hat{x}_j(k))^2\le\bigl(3(N-1)-D\bigr)(D+1)\max_{p\in\mathcal{V}}\Delta V_p(\mathbf{x}(k))$. Due to this and to $\sum_{i\in\mathcal{V}}\sum_{j\in\mathcal{V}-\mathcal{N}_i-\{i\}}b_ib_j=\sum_{i\in\mathcal{V}}\sum_{j\in\mathcal{V}}b_ib_j-\beta=N^2-\beta$, we have $\sum_{i\in\mathcal{V}}\sum_{j\in\mathcal{V}-\mathcal{N}_i-\{i\}}b_ib_j(\hat{x}_i(k)-\hat{x}_j(k))^2\le(N^2-\beta)\bigl(3(N-1)-D\bigr)(D+1)\max_{p\in\mathcal{V}}\Delta V_p(\mathbf{x}(k))$. This, along with \eqref{eq:sumsumbbxhxh2<=2aNmaxDV}, \eqref{eq:N2maxDVsumsumbbxhxh22Nsumsumbbxhxh22N}, and \eqref{eq:g=N2aN2b3N1DD12N}, implies $V(\mathbf{x}(k))\le\gamma\max_{i\in\mathcal{V}}\Delta V_i(\mathbf{x}(k))$.
\end{proof}

Because of \eqref{eq:VV=DVu}, \eqref{eq:u=argmaxDV}, and Lemma~\ref{lem:V<=gmaxDVgene}, we have $V(\mathbf{x}(k-1))-V(\mathbf{x}(k))\ge\frac{V(\mathbf{x}(k-1))}{\gamma}$ $\forall k\in\mathbb{P}$, which is exactly \eqref{eq:V<=11gV}. To prove \eqref{eq:||xx*1||<=sqrtVmax|N|211gk2} and \eqref{eq:||xhx*1||<=sqrt2Vmax|N|min|N|max|N|11gk2}, note that \eqref{eq:V<=11gV} implies $V(\mathbf{x}(k))\le(1-\frac{1}{\gamma})^kV(\mathbf{x}(0))$ $\forall k\in\mathbb{N}$. Moreover, note from \eqref{eq:V=sumcxx*2} and \eqref{eq:c=1|N|1|N|} that $V(\mathbf{x}(k))\ge(\min_{\{i,j\}\in\mathcal{E}}c_{\{i,j\}})\|\mathbf{x}(k)-x^*\mathbf{1}_L\|^2$ $\forall k\in\mathbb{N}$ where $\min_{\{i,j\}\in\mathcal{E}}c_{\{i,j\}}\ge\frac{2}{\max_{i\in\mathcal{V}}|\mathcal{N}_i|}$. Furthermore, note from \eqref{eq:V=12sumDVsumbxhx*2} and \eqref{eq:c=1|N|1|N|} that $V(\mathbf{x}(k))\ge(\min_{i\in\mathcal{V}}b_i)\|\hat{\mathbf{x}}(k)-x^*\mathbf{1}_N\|^2$ $\forall k\in\mathbb{N}$ where $\min_{i\in\mathcal{V}}b_i\ge\frac{1}{2}(1+\frac{\min_{i\in\mathcal{V}}|\mathcal{N}_i|}{\max_{i\in\mathcal{V}}|\mathcal{N}_i|})$. Thus, \eqref{eq:||xx*1||<=sqrtVmax|N|211gk2} and \eqref{eq:||xhx*1||<=sqrt2Vmax|N|min|N|max|N|11gk2} hold. To derive the bounds on $\alpha$, notice from \eqref{eq:c=1|N|1|N|} that $\frac{b_i+b_j}{c_{\{i,j\}}}=\frac{1}{2}+(1+\frac{1}{2}\sum_{\ell\in\mathcal{N}_i-\{j\}}\frac{1}{|\mathcal{N}_\ell|}+\frac{1}{2}\sum_{\ell\in\mathcal{N}_j-\{i\}}\frac{1}{|\mathcal{N}_\ell|})/(\frac{1}{|\mathcal{N}_i|}+\frac{1}{|\mathcal{N}_j|})\le\frac{1}{2}+(1+\frac{\max_{\ell\in\mathcal{V}}|\mathcal{N}_\ell|-1}{\min_{\ell\in\mathcal{V}}|\mathcal{N}_\ell|})/(\frac{2}{\max_{\ell\in\mathcal{V}}|\mathcal{N}_\ell|})\le\frac{N^2-2N+2}{2}$ $\forall\{i,j\}\in\mathcal{E}$. Similarly, it can be shown that $\frac{b_i+b_j}{c_{\{i,j\}}}\ge1$ $\forall\{i,j\}\in\mathcal{E}$. Hence, $\alpha\in[1,\frac{N^2-2N+2}{2}]$. To derive the bounds on $\beta$, observe that $\beta\le\sum_{i\in\mathcal{V}}\sum_{j\in\mathcal{V}}b_ib_j=N^2$. Also, $\sum_{i\in\mathcal{V}}\sum_{j\in\mathcal{N}_i}b_ib_j\ge2L\cdot\bigl(\frac{1}{2}(1+\frac{\min_{\ell\in\mathcal{V}}|\mathcal{N}_\ell|}{\max_{\ell\in\mathcal{V}}|\mathcal{N}_\ell|})\bigr)^2\ge\frac{L}{2}(1+\frac{1}{N-1})^2$ and $\sum_{i\in\mathcal{V}}b_i^2\ge\frac{1}{N}(\sum_{i\in\mathcal{V}}b_i)^2=N$. Therefore, $\beta\in[N+\frac{L}{2}(1+\frac{1}{N-1})^2,N^2]$. Finally, using \eqref{eq:g=N2aN2b3N1DD12N}, the bounds on $\alpha$ and $\beta$, and the properties $L\ge N-1$ and $\bigl(3(N-1)-D\bigr)(D+1)\le2N(N-1)$, we obtain $\gamma\in[\frac{N}{2}+1,N^3-2N^2+\frac{N}{2}+1]$.

\subsection{Proof of Theorem~\ref{thm:ICHAexpconvspec}}\label{ssec:proofthmICHAexpconvspec}

\begin{lemma}\label{lem:V<=gmaxDVspec}
$V(\mathbf{x}(k))\le\gamma\max_{i\in\mathcal{V}}\Delta V_i(\mathbf{x}(k))$ $\forall k\in\mathbb{N}$, where $\gamma$ is as in~\ref{enu:specpath} for a path graph with $N\ge4$, \ref{enu:speccycl} for a cycle graph, \ref{enu:specregu} for a $K$-regular graph with $K\ge2$, and~\ref{enu:specstro} for a $(N,K,\lambda,\mu)$-strongly regular graph with $\mu\ge1$.
\end{lemma}

\begin{proof}
Let $k\in\mathbb{N}$. First, suppose $\mathcal{G}$ is a path graph with $N\ge4$ and $\mathcal{E}=\{\{1,2\},\{2,3\},\ldots,\{N-1,N\}\}$. Note from \eqref{eq:x*=1Nsumy}, \eqref{eq:sumcx=sumcx}, \eqref{eq:sumcxsumc=sumyN}, and \eqref{eq:c=1|N|1|N|} that $\sum_{\{i,j\}\in\mathcal{E}}\sum_{\{p,q\}\in\mathcal{E}}c_{\{i,j\}}c_{\{p,q\}}(x_{\{i,j\}}(k)-x_{\{p,q\}}(k))^2=2N\sum_{\{i,j\}\in\mathcal{E}}c_{\{i,j\}}(x_{\{i,j\}}(k)-x^*)^2$. This, along with \eqref{eq:V=sumcxx*2} and \eqref{eq:c=1|N|1|N|}, implies that
\begin{align}
&V(\mathbf{x}(k))=\frac{1}{2N}\sum_{\{i,j\}\in\mathcal{E}}\sum_{\{p,q\}\in\mathcal{E}}c_{\{i,j\}}c_{\{p,q\}}(x_{\{i,j\}}(k)-x_{\{p,q\}}(k))^2\label{eq:V=12Nsumsumccxx2}\displaybreak[0]\\
&\quad=\frac{1}{2N}\Bigl(\sum_{\{i,j\}\in\mathcal{E}'}\sum_{\{p,q\}\in\mathcal{E}'}(x_{\{i,j\}}(k)-x_{\{p,q\}}(k))^2+3\sum_{\{i,j\}\in\mathcal{E}'}(x_{\{1,2\}}(k)-x_{\{i,j\}}(k))^2\nonumber\displaybreak[0]\\
&\quad+3\sum_{\{i,j\}\in\mathcal{E}'}(x_{\{N-1,N\}}(k)-x_{\{i,j\}}(k))^2+\frac{9}{2}(x_{\{1,2\}}(k)-x_{\{N-1,N\}}(k))^2\Bigr),\label{eq:V=12Nsumsumxx23sumxx23sumxx292xx2}
\end{align}
where $\mathcal{E}'=\mathcal{E}-\{\{1,2\},\{N-1,N\}\}$. Observe from \eqref{eq:xh=sumcxsumc}, \eqref{eq:c=1|N|1|N|}, and \eqref{eq:DV=sumcxxh2} that $(x_{\{i-1,i\}}(k)-x_{\{i,i+1\}}(k))^2=\frac{5}{3}\Delta V_i(\mathbf{x}(k))$ $\forall i\in\{2,N-1\}$ and $(x_{\{i-1,i\}}(k)-x_{\{i,i+1\}}(k))^2=2\Delta V_i(\mathbf{x}(k))$ $\forall i\in\{3,4,\ldots,N-2\}$. By the root-mean square-arithmetic mean inequality, $\sum_{\{i,j\}\in\mathcal{E}'}\sum_{\{p,q\}\in\mathcal{E}'}(x_{\{i,j\}}(k)-x_{\{p,q\}}(k))^2=2\sum_{i=2}^{N-3}\sum_{j=i+1}^{N-2}(x_{\{i,i+1\}}(k)-x_{\{j,j+1\}}(k))^2\le2\sum_{i=2}^{N-3}\sum_{j=i+1}^{N-2}(j-i)\sum_{\ell=i+1}^j(x_{\{\ell-1,\ell\}}(k)-x_{\{\ell,\ell+1\}}(k))^2=2(N-3)\sum_{i=3}^{N-2}(N-i-1)(i-2)\Delta V_i(\mathbf{x}(k))$. Moreover, $3\sum_{\{i,j\}\in\mathcal{E}'}(x_{\{1,2\}}(k)-x_{\{i,j\}}(k))^2\le3\sum_{i=2}^{N-2}(i-1)\sum_{j=2}^i(x_{\{j-1,j\}}(k)-x_{\{j,j+1\}}(k))^2=\frac{5}{2}(N-2)(N-3)\Delta V_2(\mathbf{x}(k))+3\sum_{i=3}^{N-2}(N+i-4)(N-i-1)\Delta V_i(\mathbf{x}(k))$. Similarly, $3\sum_{\{i,j\}\in\mathcal{E}'}(x_{\{N-1,N\}}(k)-x_{\{i,j\}}(k))^2\le\frac{5}{2}(N-2)(N-3)\Delta V_{N-1}(\mathbf{x}(k))+3\sum_{i=3}^{N-2}(2N-i-3)(i-2)\Delta V_i(\mathbf{x}(k))$. Finally, $\frac{9}{2}(x_{\{1,2\}}(k)-x_{\{N-1,N\}}(k))^2\le\frac{9}{2}(N-2)\sum_{i=2}^{N-1}(x_{\{i-1,i\}}(k)-x_{\{i,i+1\}}(k))^2=3(N-2)\bigl(\frac{5}{2}\Delta V_2(\mathbf{x}(k))+\frac{5}{2}\Delta V_N(\mathbf{x}(k))+3\sum_{i=3}^{N-2}\Delta V_i(\mathbf{x}(k))\bigr)$. Combining the above with \eqref{eq:V=12Nsumsumxx23sumxx23sumxx292xx2} yields $V(\mathbf{x}(k))\le\gamma\max_{i\in\mathcal{V}}\Delta V_i(\mathbf{x}(k))$ where $\gamma$ is as in~\ref{enu:specpath}.

Now suppose $\mathcal{G}$ is a cycle graph with $\mathcal{E}=\{\{1,2\},\{2,3\},\ldots,\{N-1,N\},\{N,1\}\}$. Also suppose $N$ is odd. Let $\mathbf{y}\in\mathbb{R}^N$ be a permutation of $\mathbf{x}(k)$ such that $y_{\{N,1\}}\le y_{\{1,2\}}\le y_{\{N,N-1\}}\le y_{\{2,3\}}\le y_{\{N-1,N-2\}}\le\cdots\le y_{\{\frac{N-1}{2},\frac{N+1}{2}\}}\le y_{\{\frac{N+3}{2},\frac{N+1}{2}\}}$. Then, since \eqref{eq:V=12Nsumsumccxx2} holds for any graph and due to \eqref{eq:c=1|N|1|N|}, $V(\mathbf{y})=V(\mathbf{x}(k))$. Also, due to \eqref{eq:DV=sumcxxh2} and \eqref{eq:c=1|N|1|N|}, $\max_{i\in\mathcal{V}}\Delta V_i(\mathbf{y})\le\max_{i\in\mathcal{V}}\Delta V_i(\mathbf{x}(k))$. For convenience, let $M=2\max_{i\in\mathcal{V}}\Delta V_i(\mathbf{y})$ and relabel $(y_{\{N,1\}},y_{\{1,2\}},y_{\{N,N-1\}},y_{\{2,3\}},y_{\{N-1,N-2\}},\ldots,y_{\{\frac{N-1}{2},\frac{N+1}{2}\}},\linebreak[0]y_{\{\frac{N+3}{2},\frac{N+1}{2}\}})$ as $(z_1,z_2,\ldots,z_N)$. Then, we can write $V(\mathbf{y})=\frac{1}{2N}\sum_{i=1}^N\sum_{j=1}^N(z_i-z_j)^2=\frac{1}{N}(C_1+C_2)$, where $C_1=\sum_{i=1}^{\frac{N-1}{2}}(z_1-z_{2i})^2+(z_1-z_{2i+1})^2+(z_{2i}-z_{2i+1})^2$ and $C_2=\sum_{i=1}^{\frac{N-3}{2}}\sum_{j=i+1}^{\frac{N-1}{2}}(z_{2i}-z_{2j+1})^2+(z_{2i+1}-z_{2j})^2+(z_{2i}-z_{2j})^2+(z_{2i+1}-z_{2j+1})^2$. Moreover, from \eqref{eq:xh=sumcxsumc}, \eqref{eq:c=1|N|1|N|}, and \eqref{eq:DV=sumcxxh2}, we get $z_2-z_1\le\sqrt{M}$, $z_N-z_{N-1}\le\sqrt{M}$, and $z_{i+2}-z_i\le\sqrt{M}$ $\forall i\in\{1,2,\ldots,N-2\}$. Due to the property $(a-b)^2+(a-c)^2+(b-c)^2\le2(a-c)^2$ $\forall a,b,c\in\mathbb{R}$ with $a\le b\le c$, we have $C_1\le\sum_{i=1}^{\frac{N-1}{2}}2(z_1-z_{2i+1})^2\le\sum_{i=1}^{\frac{N-1}{2}}2i^2M=\frac{N(N^2-1)}{12}M$. In addition, from the property $(a-d)^2+(b-c)^2\le(a-b)^2+(a-c)^2+(b-d)^2+(c-d)^2$ $\forall a,b,c,d\in\mathbb{R}$, we have $C_2\le \sum_{i=1}^{\frac{N-3}{2}}\sum_{j=i+1}^{\frac{N-1}{2}}2(z_{2i}-z_{2j})^2+2(z_{2i+1}-z_{2j+1})^2+(z_{2i}-z_{2i+1})^2+(z_{2j}-z_{2j+1})^2\le\sum_{i=1}^{\frac{N-3}{2}}\sum_{j=i+1}^{\frac{N-1}{2}}\bigl(4(i-j)^2M+2M\bigr)=\frac{(N-1)(N-3)(N^2+11)}{48}M$. Combining the above, we obtain $V(\mathbf{x}(k))=V(\mathbf{y})\le\frac{\gamma}{2}M\le\gamma\max_{i\in\mathcal{V}}\Delta V_i(\mathbf{x}(k))$ where $\gamma$ is as in~\ref{enu:speccycl}. Next, suppose $N$ is even. Similarly, let $\mathbf{y}\in\mathbb{R}^N$ be a permutation of $\mathbf{x}(k)$ such that $y_{\{N,1\}}\le y_{\{1,2\}}\le y_{\{N,N-1\}}\le y_{\{2,3\}}\le y_{\{N-1,N-2\}}\le\cdots\le y_{\{\frac{N}{2}-1,\frac{N}{2}\}}\le y_{\{\frac{N}{2}+2,\frac{N}{2}+1\}}\le y_{\{\frac{N}{2},\frac{N}{2}+1\}}$. Observe from \eqref{eq:V=12Nsumsumccxx2}, \eqref{eq:c=1|N|1|N|}, and \eqref{eq:DV=sumcxxh2} that $V(\mathbf{y})=V(\mathbf{x}(k))$ and $\max_{i\in\mathcal{V}}\Delta V_i(\mathbf{y})\le\max_{i\in\mathcal{V}}\Delta V_i(\mathbf{x}(k))$. As before, let $M=2\max_{i\in\mathcal{V}}\Delta V_i(\mathbf{y})$ and relabel $(y_{\{N,1\}},y_{\{1,2\}},y_{\{N,N-1\}},y_{\{2,3\}},y_{\{N-1,N-2\}},\ldots,y_{\{\frac{N}{2}-1,\frac{N}{2}\}},y_{\{\frac{N}{2}+2,\frac{N}{2}+1\}}, y_{\{\frac{N}{2},\frac{N}{2}+1\}})$ as $(z_1,z_2,\ldots,z_N)$. Then, $V(\mathbf{y})=\frac{1}{2N}\sum_{i=1}^N\sum_{j=1}^N(z_i-z_j)^2=\frac{1}{N}(C_1+C_2+C_3)$, where $C_1=\sum_{i=1}^{\frac{N}{2}-1}(z_1-z_{2i})^2+(z_1-z_{2i+1})^2+(z_{2i}-z_{2i+1})^2+(z_N-z_{2i})^2+(z_N-z_{2i+1})^2$, $C_2=\sum_{i=1}^{\frac{N}{2}-2}\sum_{j=i+1}^{\frac{N}{2}-1}(z_{2i}-z_{2j+1})^2+(z_{2i+1}-z_{2j})^2+(z_{2i}-z_{2j})^2+(z_{2i+1}-z_{2j+1})^2$, and $C_3=(z_1-z_N)^2$. Moreover, $z_2-z_1\le\sqrt{M}$, $z_N-z_{N-1}\le\sqrt{M}$, and $z_{i+2}-z_i\le\sqrt{M}$ $\forall i\in\{1,2,\ldots,N-2\}$. Using the above properties, it can be shown that $C_1\le C_1+\sum_{i=1}^{\frac{N}{2}-1}(z_{2i}-z_{2i+1})^2\le\sum_{i=1}^{\frac{N}{2}-1}2(z_1-z_{2i+1})^2+2(z_N-z_{2i})^2\le\sum_{i=1}^{\frac{N}{2}-1}2i^2M+2(\frac{N}{2}-i)^2M=\frac{N(N-1)(N-2)}{6}M$, $C_2\le\sum_{i=1}^{\frac{N}{2}-2}\sum_{j=i+1}^{\frac{N}{2}-1}\bigl(4(i-j)^2M+2M\bigr)=\frac{(N-2)(N-4)(N^2-2N+12)}{48}M$, and $C_3\le\frac{N^2}{4}M$. It follows that $V(\mathbf{x}(k))\le\gamma\max_{i\in\mathcal{V}}\Delta V_i(\mathbf{x}(k))$ where $\gamma$ is as in~\ref{enu:speccycl}.

Next, suppose $\mathcal{G}$ is a $K$-regular graph with $K\ge2$. Due to \eqref{eq:c=1|N|1|N|} and \eqref{eq:DV=sumcxxh2}, $\sum_{i\in\mathcal{V}}\Delta V_i(\mathbf{x}(k))=\frac{2}{K}\sum_{\{i,j\}\in\mathcal{E}}(\hat{x}_i(k)-x_{\{i,j\}}(k))^2+(\hat{x}_j(k)-x_{\{i,j\}}(k))^2\ge\frac{1}{K}\sum_{\{i,j\}\in\mathcal{E}}(\hat{x}_i(k)-\hat{x}_j(k))^2$, implying that
\begin{align}
\sum_{i\in\mathcal{V}}\sum_{j\in\mathcal{N}_i}(\hat{x}_i(k)-\hat{x}_j(k))^2\le2K\sum_{i\in\mathcal{V}}\Delta V_i(\mathbf{x}(k)).\label{eq:sumsumxhxh2<=2KsumDV}
\end{align}
Again, because of \eqref{eq:c=1|N|1|N|} and \eqref{eq:DV=sumcxxh2}, $\forall i\in\mathcal{V}$, $\forall j\in\mathcal{N}_i$, $(x_{\{i,j\}}(k)-\hat{x}_i(k))^2\le\frac{K}{2}\max_{p\in\mathcal{V}}\Delta V_p(\mathbf{x}(k))$. Moreover, $\forall i\in\mathcal{V}$, $\forall j,\ell\in\mathcal{N}_i$ with $j\neq\ell$, $(x_{\{i,j\}}(k)-x_{\{i,\ell\}}(k))^2\le2\bigl((x_{\{i,j\}}(k)-\hat{x}_i(k))^2+(x_{\{i,\ell\}}(k)-\hat{x}_i(k))^2\bigr)\le K\max_{p\in\mathcal{V}}\Delta V_p(\mathbf{x}(k))$. Via the preceding two inequalities and the root-mean square-arithmetic mean inequality, it can be shown that $\forall i\in\mathcal{V}$, $\forall j\in\mathcal{V}-\mathcal{N}_i-\{i\}$, $(\hat{x}_i(k)-\hat{x}_j(k))^2\le(D+1)(\frac{K}{2}\max_{p\in\mathcal{V}}\Delta V_p(\mathbf{x}(k))\cdot2+K\max_{p\in\mathcal{V}}\Delta V_p(\mathbf{x}(k))\cdot(D-1))=KD(D+1)\max_{p\in\mathcal{V}}\Delta V_p(\mathbf{x}(k))$. It then follows from \eqref{eq:N2maxDVsumsumbbxhxh22Nsumsumbbxhxh22N}, \eqref{eq:c=1|N|1|N|}, and \eqref{eq:sumsumxhxh2<=2KsumDV} that $V(\mathbf{x}(k))\le\gamma\max_{i\in\mathcal{V}}\Delta V_i(\mathbf{x}(k))$ where $\gamma$ is as in~\ref{enu:specregu}.

Finally, suppose $\mathcal{G}$ is a $(N,K,\lambda,\mu)$-strongly regular graph with $\mu\ge1$, which means that it is a $K$-regular graph with $K\ge2$ and with every two non-adjacent nodes having $\mu$ common neighbors. For every $i\in\mathcal{V}$ and $j\in\mathcal{V}-\mathcal{N}_i-\{i\}$, let $\{q_{ij1},q_{ij2},\ldots,q_{ij\mu}\}=\mathcal{N}_i\cap\mathcal{N}_j$. Then, from \eqref{eq:c=1|N|1|N|} and \eqref{eq:DV=sumcxxh2}, $\mu\sum_{i\in\mathcal{V}}\sum_{j\in\mathcal{V}-\mathcal{N}_i-\{i\}}(\hat{x}_i(k)-\hat{x}_j(k))^2=\sum_{i\in\mathcal{V}}\sum_{j\in\mathcal{V}-\mathcal{N}_i-\{i\}}\sum_{\ell=1}^\mu(\hat{x}_i(k)-\hat{x}_j(k))^2\le4\sum_{i\in\mathcal{V}}\sum_{j\in\mathcal{V}-\mathcal{N}_i-\{i\}}\sum_{\ell=1}^\mu\bigl((\hat{x}_i(k)-x_{\{i,q_{ij\ell}\}}(k))^2+(x_{\{i,q_{ij\ell}\}}(k)-\hat{x}_{q_{ij\ell}}(k))^2+(\hat{x}_{q_{ij\ell}}(k)-x_{\{j,q_{ij\ell}\}}(k))^2+(x_{\{j,q_{ij\ell}\}}(k)-\hat{x}_j(k))^2\bigr)\le2K\sum_{i\in\mathcal{V}}\sum_{j\in\mathcal{V}-\mathcal{N}_i-\{i\}}\bigl(\Delta V_i(\mathbf{x}(k))+\sum_{\ell=1}^\mu\Delta V_{q_{ij\ell}}(\mathbf{x}(k))+\Delta V_j(\mathbf{x}(k))\bigr)\le2KN(N-K-1)(2+\mu)\max_{p\in\mathcal{V}}\Delta V_p(\mathbf{x}(k))$. This, along with \eqref{eq:N2maxDVsumsumbbxhxh22Nsumsumbbxhxh22N}, \eqref{eq:c=1|N|1|N|}, and \eqref{eq:sumsumxhxh2<=2KsumDV}, implies that $V(\mathbf{x}(k))\le\gamma\max_{i\in\mathcal{V}}\Delta V_i(\mathbf{x}(k))$ where $\gamma$ is as in~\ref{enu:specstro}.
\end{proof}

Note that in the proof of Theorem~\ref{thm:ICHAexpconvgene} in Appendix~\ref{ssec:proofthmICHAexpconvgene}, Lemma~\ref{lem:V<=gmaxDVgene} is used to derive \eqref{eq:V<=11gV}--\eqref{eq:||xhx*1||<=sqrt2Vmax|N|min|N|max|N|11gk2}. In the same way, \eqref{eq:V<=11gV}--\eqref{eq:||xhx*1||<=sqrt2Vmax|N|min|N|max|N|11gk2} can be derived using Lemma~\ref{lem:V<=gmaxDVspec}, completing the proof of Theorem~\ref{thm:ICHAexpconvspec}.

\subsection{Proof of Theorem~\ref{thm:CHAexpconv}}\label{ssec:proofthmCHAexpconv}

Let $\gamma$ be as in \eqref{eq:g=N2aN2b3N1DD12N} for a general graph or as in~\ref{enu:specpath}--\ref{enu:specstro} for a specific graph. Note that Lemmas~\ref{lem:V<=gmaxDVgene} and~\ref{lem:V<=gmaxDVspec} are independent of $(u(k))_{k=1}^\infty$ and, thus, hold for CHA as well. Hence,
\begin{align}
V(\mathbf{x}(k))\le\gamma\max_{i\in\mathcal{V}}\Delta V_i(\mathbf{x}(k)),\quad\forall k\in\mathbb{N}.\label{eq:V<=gmaxDV}
\end{align}
Next, analyzing Algorithm~\ref{alg:CHA} with $\varepsilon(v)=0$ $\forall v\in[0,\infty)$, we see that
\begin{align}
&V(\mathbf{x}(1))=V(\mathbf{x}(0))-\max_{i\in\mathcal{V}}\Delta V_i(\mathbf{x}(0)),\label{eq:V=VmaxDV}\displaybreak[0]\\
&V(\mathbf{x}(k+1))\le V(\mathbf{x}(k))-\min\{\max_{i\in\mathcal{V}}\Delta V_i(\mathbf{x}(k)),\Phi^{-1}(t(k))\},\quad\forall k\in\mathbb{P},\label{eq:V<=VminmaxDVPinvt}\displaybreak[0]\\
&t(k+1)=\max\{\Phi(\max_{i\in\mathcal{V}}\Delta V_i(\mathbf{x}(k))),t(k)\},\quad\forall k\in\mathbb{N}.\label{eq:t=maxPmaxDVt}
\end{align}
With \eqref{eq:V<=gmaxDV}--\eqref{eq:t=maxPmaxDVt}, we now show by induction that $\forall k\in\mathbb{P}$, $V(\mathbf{x}(k))\le(1-\frac{1}{\gamma})V(\mathbf{x}(k-1))$ and $t(k)\le\Phi(\frac{V(\mathbf{x}(k-1))}{\gamma})$. Let $k=1$. Then, because of \eqref{eq:V<=gmaxDV}, \eqref{eq:V=VmaxDV}, and \eqref{eq:t=maxPmaxDVt} and because $\Phi$ is strictly decreasing, we have $V(\mathbf{x}(1))\le(1-\frac{1}{\gamma})V(\mathbf{x}(0))$ and $t(1)=\Phi(\max_{i\in\mathcal{V}}\Delta V_i(\mathbf{x}(0)))\le\Phi(\frac{V(\mathbf{x}(0))}{\gamma})$. Next, let $k\ge1$ and suppose $V(\mathbf{x}(k))\le(1-\frac{1}{\gamma})V(\mathbf{x}(k-1))$ and $t(k)\le\Phi(\frac{V(\mathbf{x}(k-1))}{\gamma})$. To show that $V(\mathbf{x}(k+1))\le(1-\frac{1}{\gamma})V(\mathbf{x}(k))$ and $t(k+1)\le\Phi(\frac{V(\mathbf{x}(k))}{\gamma})$, consider the following two cases: (i) $\max_{i\in\mathcal{V}}\Delta V_i(\mathbf{x}(k))<\Phi^{-1}(t(k))$ and (ii) $\max_{i\in\mathcal{V}}\Delta V_i(\mathbf{x}(k))\ge\Phi^{-1}(t(k))$. For case~(i), due to \eqref{eq:V<=gmaxDV}, \eqref{eq:V<=VminmaxDVPinvt}, and \eqref{eq:t=maxPmaxDVt}, we have $V(\mathbf{x}(k+1))\le V(\mathbf{x}(k))-\max_{i\in\mathcal{V}}\Delta V_i(\mathbf{x}(k))\le(1-\frac{1}{\gamma})V(\mathbf{x}(k))$ and $t(k+1)=\Phi(\max_{i\in\mathcal{V}}\Delta V_i(\mathbf{x}(k)))\le\Phi(\frac{V(\mathbf{x}(k))}{\gamma})$. For case~(ii), due to \eqref{eq:V<=VminmaxDVPinvt}, \eqref{eq:t=maxPmaxDVt}, and the hypothesis, we have $V(\mathbf{x}(k+1))\le V(\mathbf{x}(k))-\Phi^{-1}(t(k))\le V(\mathbf{x}(k))-\frac{V(\mathbf{x}(k-1))}{\gamma}\le V(\mathbf{x}(k))-\frac{V(\mathbf{x}(k))}{\gamma(1-\frac{1}{\gamma})}\le(1-\frac{1}{\gamma})V(\mathbf{x}(k))$ and $t(k+1)=t(k)\le\Phi(\frac{V(\mathbf{x}(k-1))}{\gamma})\le\Phi(\frac{V(\mathbf{x}(k))}{\gamma(1-\frac{1}{\gamma})})\le\Phi(\frac{V(\mathbf{x}(k))}{\gamma})$. This completes the proof by induction. It follows that \eqref{eq:V<=11gV} and therefore \eqref{eq:||xx*1||<=sqrtVmax|N|211gk2} and \eqref{eq:||xhx*1||<=sqrt2Vmax|N|min|N|max|N|11gk2} hold, so that Theorems~\ref{thm:ICHAexpconvgene} and~\ref{thm:ICHAexpconvspec} hold verbatim here. Next, observe from \eqref{eq:t=maxPmaxDVt} that $(t(k))_{k=0}^\infty$ is non-decreasing. To show that $\lim_{k\rightarrow\infty}t(k)=\infty$, assume to the contrary that $\exists \bar{t}\in(0,\infty)$ such that $t(k)\le\bar{t}$ $\forall k\in\mathbb{N}$. For each $k\in\mathbb{P}$, reconsider the above two cases. Because of \eqref{eq:V<=VminmaxDVPinvt} and \eqref{eq:t=maxPmaxDVt}, for case~(i), $V(\mathbf{x}(k))-V(\mathbf{x}(k+1))\ge\max_{i\in\mathcal{V}}\Delta V_i(\mathbf{x}(k))=\Phi^{-1}(t(k+1))\ge\Phi^{-1}(\bar{t})$. Similarly, for case~(ii), $V(\mathbf{x}(k))-V(\mathbf{x}(k+1))\ge\Phi^{-1}(t(k))\ge\Phi^{-1}(\bar{t})$. Combining these two cases, we get $V(\mathbf{x}(k+1))\le V(\mathbf{x}(1))-k\Phi^{-1}(\bar{t})$ $\forall k\in\mathbb{N}$. Since $\Phi^{-1}(\bar{t})>0$, $V(\mathbf{x}(k+1))<0$ for sufficiently large $k$, which is a contradiction. Thus, $\lim_{k\rightarrow\infty}t(k)=\infty$. Finally, from the statement shown earlier by induction, we obtain $V(\mathbf{x}(k))\le(1-\frac{1}{\gamma})\cdot\gamma\Phi^{-1}(t(k))=(\gamma-1)\Phi^{-1}(t(k))$ $\forall k\in\mathbb{P}$.

\bibliographystyle{IEEEtran}
\bibliography{paper}

\end{document}